\documentclass[12pt,amstex]{amsart}
\usepackage{amsmath}
\usepackage{amscd}
\usepackage{amsfonts}
\usepackage{amssymb}
\usepackage{amsthm}

  =4.mm
\newcommand{\rar}{\rightarrow}
\newcommand{\calg}{\mathcal}
\newcommand{\hb}{\widehat{B}}
\newcommand{\hc}{\widehat{C}}
\newcommand{\ilim}{\text{lim } }
\newcommand{\rhh}{\text{RHom}}
\newcommand{\rhhh}{\text{R}\mathcal{H}om}

\newcommand{\chh}{{\mathcal{H}om}^{\text{cont}}}
\newcommand{\dd}{\text{D}_{poly}}
\newcommand{\strc}{\mathcal{O}_{X}}

\newcommand{\tenp}{\otimes_{\mathcal{O}_{X}}}

\newcommand{\dcat}{\text{D}^{+}(X)}

\newcommand{\dcatb}{\text{D}^{+}(X)}
\newcommand{\bbcc}{\text{Ch}^{+}({{\mathcal O}_X}-{\text{mod}})}

\newcommand{\exal}[2]{ \wedge^{#2} \Omega_{#1} [#2]}
\newcommand{\texal}[2]{ \wedge^{#2} T_{#1} [-#2]}
\newcommand{\hkr}{\text{I}_{HKR}}

\newcommand{\pbw}{\text{I}_{PBW}}
\newcommand{\bhkr}{\text{J}}

\newcommand{\sss}{\text{Sym}}
\newcommand{\ttt}{\text{t}}
\newcommand{\enn}{\text{End}}
\newcommand{\ennnn}{{\mathcal E}\text{nd}}
\newcommand{\tot}{\text{Tot}}
\newcommand{\tdd}{\triangledown}
\newcommand{\chn}{\text{ch}}
\newcommand{\caph}{\text{H}}
\newcommand{\exx}{\text{exp}}
\newcommand{\adj}{\text{ad}}
\newcommand{\id}{\text{id}}
\newcommand{\pll}{\Pi \text{LIE}}
\newcommand{\stdp}{\text{END}_{T[-1]}}

\newcommand{\frel}[1]{{\mathcal L}(#1)}
\newcommand{\ufrel}[1]{\mathcal{T}(#1)}
\newcommand{\smap}[1]{\text{J}(#1)}
\newcommand{\bsmap}[1]{\text{B}(#1)}

\newtheorem{thm}{Theorem}
\newtheorem*{thm1}{Theorem 1}
\newtheorem*{thm2}{Theorem 2}
\newtheorem*{thm3}{Theorem 3}
\newtheorem*{thm4}{Theorem 4}
\newtheorem*{thm5}{Theorem 5}

\newtheorem{prop}{Proposition}
\newtheorem{cor}{Corollary}
\newtheorem{lem}{Lemma}

\title{ The big Chern classes and the Chern character}
\author{Ajay C. Ramadoss}

 \address{Department of Mathematics, University of Oklahoma \newline
 \text{        }       621 Elm Avenue, Norman,  OK-73019}

  \email{aramadoss@math.ou.edu}
\begin{document}

\maketitle
\begin{abstract}
Let $X$ be a smooth scheme over a field of characteristic $0$. Let
$\dd^{\bullet}(X)$ be the complex of polydifferential operators on
$X$ equipped with Hochschild co-boundary. Let $L(\dd^1(X))$ be the
free Lie algebra generated over $\strc$ by $\dd^1(X)$ concentrated
in degree $1$ equipped with Hochschild co-boundary. We have a
symmetrization map $I: \oplus_k \sss^k(L(\dd^1(X))) \rar
\dd^{\bullet}(X)$. Theorem 1 of this
paper measures how the map $I$ fails to commute with multiplication. \\

We recall from Kapranov [6] that the Atiyah class of the tangent
bundle yields a map $\alpha_{T_X}: T_X[-1] \otimes T_X[-1] \rar
T_X[-1]$ in $\dcat$ making $T_X[-1]$ a Lie algebra object in
$\dcat$. Theorem 2 enables us to realize the Atiyah class of $T_X$
as a honest map of complexes. Theorem 2 of this paper says that
$T_X[-1]$ is quasi-isomorphic to $L(\dd^1(X))$. It further states
that the natural Lie bracket on $L(\dd^1(X))$ represents the Atiyah
class of $T_X$. An immediate consequence of Theorems 1 and 2 is
Corollary 1, a result "dual" to Theorem 1 of Markarian [3] that
measures how the Hochschild-Kostant-Rosenberg quasi-isomorphism fails to commute with multiplication.\\

In order to understand Theorem 1 conceptually, we prove a theorem (Theorem 3) stating that $\dd^{\bullet}(X)$ is the universal enveloping
algebra of $T_X[-1]$ in $\dcat$. At this juncture, we recall from Kapranov [6] that if $E$ is a vector bundle of $E$, the Atiyah class of $E$
equips $E$ with the structure of a module over the Lie algebra $T_X[-1]$ in $\dcat$. An easy consequence of Theorem 3 is Theorem 4, which
interprets the
Chern character $E$ as the "character of the representation $E$ of $T_X[-1]$" and gives a description of the big Chern classes of $E$. \\

Finally, Theorem 4 along with Theorem 1 is used to give an explicit formula (Theorem 5) expressing the big Chern classes of $E$ in terms of the
components of the Chern character of $E$. \\

{\it Keywords :} Complex of polydifferential operators, Hochschild-Kostant-Rosenberg quasi-isomorphism,
 Atiyah class, big Chern classes, Chern character, derived category, lie algebra, universal enveloping algebra . \\

Mathematics subject classification 2000: 14F99,16E40,18G99.

\end{abstract}

\section{Introduction}

\subsection{Outline of main results}

This paper is a result of an effort to understand the works of Markarian [3] and Caldararu [10]. Another goal was to see whether the works of
Markarian [3] and Caldararu [10] led to an explicit formula relating the big Chern classes to the Chern character. \\

We begin by outlining the main results in this paper. Let $X$ be a smooth scheme over a field of characteristic $0$. Let $\dd^{\bullet}(X)$
denote the complex of poly-differential operators on $X$, with Hochschild co-boundary. We denote the sheaf of differential operators on $X$ by
$\dd^1(X)$. Let $L(\dd^1(X))$ denote the free Lie algebra generated over $\strc$ by $\dd^1(X)$ concentrated in degree $1$ , equipped with the
Hochschild co-boundary. There is a symmetrization map $$I:\oplus_k \sss^k(L(\dd^1(X))) \rar \dd^{\bullet}(X)$$. $I$ is an isomorphism of
complexes of $\strc$-modules. Let \\ $m:\dd^{\bullet}(X) \otimes \dd^{\bullet}(X) \rar \dd^{\bullet}(X)$ denote the multiplication on
$\dd^{\bullet}(X)$ and let $\mu$ denote the natural product on $\oplus_k \sss^k(L(\dd^1(X)))$. Let $$\adj: \oplus_k \sss^k(L(\dd^1(X))) \otimes
L(\dd^1(X)) \rar \oplus_k \sss^k(L(\dd^1(X)))$$ denote the right adjoint action of $L(\dd^1(X))$ on $\oplus_k \sss^k(L(\dd^1(X)))$. In Section 5
,we describe a map $$\omega: \oplus_k \sss^k(L(\dd^1(X))) \otimes L(\dd^1(X)) \rar \oplus_k \sss^k(L(\dd^1(X))) \otimes L(\dd^1(X))$$ such that
$\mu \circ \omega =
\adj$. Let $\bbcc$ denote the category of bounded below complexes of $\strc$-modules. We then have the following Theorem \\

\begin{thm}
The following diagram commutes in $\bbcc$ \\

$$\begin{CD}
\dd^{\bullet}(X) \otimes L(\dd^1(X)) @> m >> \dd^{\bullet}(X) \\
@AA{I \otimes id}A           @A{I}AA \\
\sss^{\bullet}(L(\dd^1(X))) \otimes L(\dd^1(X)) @> \mu \circ \frac{\omega}{1-e^{-\omega}}>> \sss^{\bullet}(L(\dd^1(X)))\\
\end{CD} $$

\end{thm}

The next sections of this paper are devoted to a conceptual understanding of the above theorem and its corollaries. Let $\dcat$ denote the
derived category of bounded below complexes of $\strc$-modules with coherent cohomology. Recall that the Atiyah class of a vector bundle $E$ is
an element in $\text{Hom}_{\dcat}(E \otimes T_X[-1], E)$. In particular, the Atiyah class $\alpha_{TX}$ of the tangent bundle of $X$ is an
element in $$ \text{Hom}_{\dcat}(T_X \otimes T_X[-1], T_X) = \text{Hom}_{\dcat}(T_X[-1] \otimes T_X[-1] , T_X[-1])$$. Let $\hkr$ denote the
Hochschild-Kostant-Rosenberg quasi-isomorphism
$$\hkr: \oplus_i \texal{X}{i} \rar \dd^{\bullet}(X)$$ In Section 4.2, we define a map of complexes $\beta:T_X[-1] \rar L(\dd^1(X))$ such that $$
I \circ \sss^{\bullet} \beta = \hkr $$ and show that
$\beta: T_X[-1] \rar L(\dd^1(X))$ is a quasi-isomorphism.  Moreover as proven in Section 6, \\

\begin{thm}

The following diagram commutes in $\dcat$ \\

$$\begin{CD} T_X[-1] \tenp T_X[-1] @> \beta \otimes \beta >> L(\dd^1(X)) \tenp
  L(\dd^1(X)) \\
  @V \alpha(T_X) VV                    @VV [,] V \\
  T_X[-1] @> \beta >> L(\dd^1(X)) \\
\end{CD}$$
\\

\end{thm}

Note that all arrows in the diagram above except $\alpha_{T_X}$ are
honest maps in  \\ $\bbcc$. This Theorem says that the natural Lie
bracket on $L(\dd^1(X))$ realizes the Atiyah class of $T_X$ as a map
in $\bbcc$. It follows from the fact that $\beta$ is a
quasi-isomorphism that the map $\omega$ yields a map
$$\bar{\omega}: \oplus_i \texal{X}{i} \otimes T_X[-1]  \rar \oplus_i \texal{X}{i} \otimes T_X[-1]$$ in $\dcat$ . By an abuse of notation, let $\mu$ also
denote the natural product in $\oplus_i \texal{X}{i}$. An immediate consequence
of Theorems 1 and 2 is the following corollary \\

\begin{cor}
The following diagram commutes in $\dcat$ \\
$$\begin{CD} \dd^{\bullet}(X) \otimes \dd^{\bullet}(X) @> m >> \dd^{\bullet}(X) \\
@AA \hkr \otimes \hkr A            @A \hkr AA \\
\oplus_i \texal{X}{i} \otimes T_X[-1] @> \mu \circ
\frac{\bar{\omega}}{1-e^{-\bar{\omega}}} >> \oplus_i \texal{X}{i} \\
\end{CD}$$ \\
\end{cor}

This is a result "dual" to Theorem 1 of [3]. In another paper [16]
we use this to prove the relative Riemann-Roch theorem, thereby
completely explaining [3]. One can think of the
Hochschild-Kostant-Rosenberg map as a symmetrization map from
$\oplus_i \texal{X}{i}$ to $\dd^{\bullet}(X)$. Corollary 1 thus
tells us that the error term measuring how this map fails to commute
with multiplication is "$d(\text{exp}^{-1}$-like". In order to
understand the classical situation of which Corollary 1 is an
analog, we have the following theorem.

\begin{thm}
$\dd^{\bullet}(X)$ is the universal enveloping algebra of $T_X[-1]$ in $\dcat$. In other words, let $A$ be an associative algebra in $\dcat$.
If $f: T_X[-1] \rar A$ is a morphism in $\dcat$ making the following diagram commute in $\dcat$ \\

 $$\begin{CD}
 T_X[-1] \tenp T_X[-1] @> \alpha(T_X) >> T_X[-1] \\
 @V f \otimes f VV                     @VVfV \\
 A \tenp A @> [, ] >> A \\
 \end{CD}$$ \\

 then $\exists$ a unique morphism  $\bar{f}: \dd^{\bullet}(X) \rar A$ of algebras in
 $\dcat$ so that the composite $\bar{f} \circ \hkr = f$.
 \end{thm}

This makes the parallel between Theorem 1, Corollary 1 and their
classical analogs more explicit. In fact, as explained in Section 5,
the classical analog of Theorem 1 is a commutative diagram
equivalent to the formula $d({\text{exp}}^{-1}) =
\frac{\adj}{1-e^{-\adj}}$. Theorem 3
shows us that Corollary 1 is also analogous to the same classical result.  \\

In the original version of this paper, Theorem 3 was proven first. Corollary 1 was then interpreted as the analog for the Lie algebra $T_X[-1]$
of its classical version, which followed from the formula $d({\text{exp}}^{-1}) = \frac{\adj}{1-e^{-\adj}}$. However, no explicit details of its
proof were given. The same interpretation of an equivalent result has been offered by Markarian [14] in an Arxiv preprint subsequent to the
first Arxiv version of this paper. Again, as in the original version of this paper, hardly any further details were offered. Theorem 1 removes
these shortcomings. It is also stronger. \\

Theorem 3 is equivalent to a result claimed by Roberts [7]. It was however, not proven in printed or online literature available to the author
when the first version of this paper was written. In a paper [11] that appeared on the Arxiv after the first version of this paper was uploaded
on the Arxiv, Roberts and Willerton prove an equivalent result. Their proof is however, very different from the proof here. The proof here is
more explicit in the sense that $\bar{f}$ as in Theorem 3 is directly constructed in our proof.  \\

There are other, more serious applications of Theorem 3. If $E$ is a
vector bundle on $X$, the Atiyah class of $E$, $\alpha(E)$ is an
element of \\ $\text{Hom}_{\dcat}(E \otimes T_X[-1] , E) =
\text{Hom}_{\dcat}(T_X[-1], \ennnn(E))$. By [6], $\alpha(E)$ equips
$E$ with the structure of a module over the Lie algebra $T_X[-1]$ in
$\dcat$. Thus, \\ $\alpha(E): T_X[-1] \rar \ennnn(E)$ is a morphism
of Lie algebras in $\dcat$. By Theorem 3,
 there exists a morphism $\theta_E:\dd^{\bullet}(X) \rar \ennnn(E)$ of algebras in $\dcat$ lifting $\alpha(E)$. We also have a map
 $tr:\ennnn(E) \rar \strc$. Let $\varphi_E = tr \circ \theta_E$. Let $p:\oplus_i \texal{X}{i} \rar \oplus_n T_X^{\otimes n}[-n]$ be the symmetrization map.
 There is a map of complexes $\bhkr: \oplus_n T_X^{\otimes n}[-n] \rar \dd^{\bullet}(X)$, such that
such that $\hkr = \bhkr \circ p$. Let $\tilde{\ttt_k}(E)$ denote $\alpha(E)^{\circ k} \in \text{Hom}_{\dcat}(T_X^{\otimes k}[-k] , \ennnn(E))$.
Let $\tilde{\chn_k(E)} = \tilde{\ttt_k(E)} \circ p $. Let $\ttt_k(E)$ denote the $k$ th big Chern class of $E$ and let $\chn_k(E)$ denote the
$k$ th component of the Chern character of $E$. The following easy consequence of Theorem 3 is stated as a theorem in its own right.

\begin{thm}

(i) $$\tilde{\ttt_k(E)} = \theta_E \circ \bhkr $$ (ii) $$ \tilde{\chn_k(E)} = \theta_E \circ \hkr  $$ (iii) $$ \oplus_n \ttt_n(E) = \varphi_E
\circ \bhkr $$ (iv) $$ \chn(E) = \varphi_E \circ \hkr $$

\end{thm}

Part (iv) of Theorem 4 interprets the Chern character of $E$ as the
"character of the representation $E$ of $T_X[-1]$". We will comment
on this aspect in greater detail in Section 8. Theorem 4 is similar
to Theorem 4.5 of Caldararu [10]. However, that result does not lend
itself to our interpretation of the Chern character as directly as
Theorem 4 does. Further, even more interesting applications of
Theorem 3 may be found in [7], but they are beyond the scope of this paper. \\

Note that we have a PROP $\stdp$ such that  \\ $\stdp(n,m) =
\text{Hom}_{\dcat}(T_X^{\otimes n}[-n], T_X^{\otimes m}[-m])$. Let
$\Psi \in \oplus_{m \leq n} \stdp(n,m)$ be the element of $\stdp$
given in Section 9. Let $\Psi_{kl}$ denote the component of $\Psi$
in $\stdp(k,l)$. Let $\pi: \oplus_n T_X^{\otimes n}[-n] \rar
\oplus_i \texal{X}{i}$ be the standard projection. Then,

\begin{thm}

(i) $$ \tilde{\ttt_k}(E) =\tilde{\chn_k}(E) \circ \pi + \sum_{l<k}    \tilde{\chn_l(E)} \circ \pi \circ \Psi_{kl} $$ (ii) $$\ttt_k(E) =
\chn_k(E) \circ \pi + \sum_{l < k} \chn_l(E) \circ \pi \circ \Psi_{kl}$$
\\

\end{thm}

We remark here that $\Psi$ has been described by an explicit, albeit lengthy formula in Section 9 . This expresses the big Chern classes of a
vector bundle on an arbitrary smooth scheme over a field of characteristic $0$ in terms of the components of its Chern character. The existence
of a formula similar to (ii) was proven in the author's thesis for vector bundles on projective varieties over a field of characteristic $0$.
The proof there was entirely different. It crucially required the existence of an ample line bundle on $X$. The proof there therefore did not
generalize to arbitrary smooth schemes ,smooth complex manifolds etc unlike the proof here. Further, it was difficult to see the formula arising
out of the trace applied to an almost identical formula akin to (i) of Theorem 5. Also, the explicit description of $\Psi$ was not given in the
author's thesis [8].

\subsection{Structure of this paper}

Sections 2 and 3 are introductory and describe the basic properties of $\dd^{\bullet}(X)$. In particular, Section 3 describes the Hopf -algebra
structure of $\dd^{\bullet}(X)$ in $\bbcc$. \\

Section 4 proves a key lemma (Lemma 1) stating that the
symmetrization map \\ $I:\oplus_k \sss^k(L(\dd^1(X))) \rar
\dd^{\bullet}(X)$ is an isomorphism in $\bbcc$. This is done by
first showing that it is a map of complexes of $\strc$-modules,
followed by showing that it is an
isomorphism of graded $\strc$-modules. \\

Section 4 also recalls the definition of the Hochschild-Kostant-Rosenberg (HKR) quasi-isomorphism $\hkr$ from Yekutieli and shows that there is
a
quasi-isomorphism \\ $\beta: T_X[-1] \rar L(\dd^1(X))$ so that $I \circ \sss^{\bullet} \beta =\hkr$. This is later used in Section 6. \\

Section 5 states and proves Theorem 1. Once again, this is done in stages. The first stage involves an explicit calculation showing that all
maps involved commute with the relevant differentials. This enables us to "forget" the differentials in the complexes involved. We then only
need to show that the diagram involved commutes in the category of graded $\strc$-modules. \\

Let $V$ be a vector space over a field $K$ of characteristic $0$. Let $L(V)$ be the (graded) free Lie algebra generated over $K$ by $V$
concentrated in degree $1$. The universal enveloping algebra of $L(V)$ is the tensor algebra $T(V)$ of $V$. The second stage in proving Theorem
1 involves reducing the problem further to the problem of finding the error term that measures how the PBW map from $\sss^{\bullet}(L(V))$ to
$T(V)$ fails to commute with multiplication. This is carried out in the appendix at the end of this paper. \\

Section 5 also contains a lemma (Lemma 2) that furthers our understanding of Theorem 1. This states that $\dd^{\bullet}(X)$ is the universal
enveloping algebra of $L(\dd^1(X))$ in $\bbcc$. This enables us to interpret $I$ as a PBW-map from $\sss^{\bullet}(L(\dd^1(X)))$ to $\dd^{\bullet}(X)$.
Theorem 1 then says that the error term that measures the failure of $I$ to commute with multiplication is "$d(\text{exp}^{-1})$ like".\\

Section 6 recalls the definition of the Atiyah class of a perfect complex of $\strc$-modules on $X$. Theorem 2, which states that the Lie
bracket on $L(\dd^1(X))$ realizes the Atiyah class of $T_X$ as a map of complexes of $\strc$-modules, is also proven here by an explicit
computation. Theorem 1 and Theorem 2 immediately imply Corollary 1, which is also stated in Section 6. \\

Section 7 is devoted to the proof of Theorem 3. We have attempted a careful and self-contained treatment of Theorem 3 in this section.
Theorem 3 enables us to interpret the HKR-quasi-isomorphism as a PBW-map. Corollary 1 then says that the error term that measures the failure
of the HKR map to commute with multiplication in $\dcat$ is "$d(\text{exp}^{-1})$ like".   \\

Section 8 is used to state and prove Theorem 4. A spinoff of this result is a new conceptual proof of  result already proven in the author's
thesis (Corollary 6). This result states that the big Chern classes commute with Adams operations. On the other hand, one had the representation
theoretic identity $\chi_{\psi^p E}(g) = \chi_E(g^p)$ for any element $g$ of a group $G$ and for any representation $E$ of $G$. The parallel
between the fact that the big Chern classes commute with Adams operations and the identity $\chi_{\psi^p E}(g) = \chi_E(g^p)$ is
made transparent by this proof. \\

Section 9 is devoted to describing a formula for the element $\Psi$ of the PROP $\stdp$ mentioned before stating Theorem 5, and then proving
Theorem 5. \\

There is an appendix at the end of this paper. Let $K$ be a field of
characteristic $0$. Let $V$ be a graded vector space concentrated in
degree $1$.  Let $L(V)$ be the (graded) free Lie algebra generated
by $V$ over $K$. Let $T(V)$ be the tensor algebra of $V$. The
appendix devoted to a result (Theorem 6) about the error term
measuring how the PBW map from $\sss^{\bullet}(L(V))$ to $T(V)$
fails to commute with multiplication. It says that this error term
is "$d(\text{exp}^{-1})$like". The precise statement and proof are
in the Appendix. Even though we expect this to be standard, such a
result is standard only for the case when $V$ is concentrated in
degree $0$ as far as I know. This is the reason for this result to
be
included as a theorem in the Appendix to this paper. \\

\subsection*{Notation}

Throughout this paper , $\bbcc$ will denote the category of bounded below complexes of $\strc$-modules. $\dcat$ will denote the bounded below
derived category of complexes of $\strc$-modules with coherent co-homology. \\

\section{The complete Hochschild chain and cochain complexes}

The purpose of this section is to recall definitions of and facts about the completed Hochschild chain and cochain complexes. Most of the
material in this section is recalled from Yekutieli's paper [1]. The notation also follows the same source closely.
\\

Throughout this article, we shall work with smooth schemes over fields of characteristic $0$. Let $X$ be a smooth separated scheme over a field
$K$ of characteristic $0$. We have the (closed) diagonal embedding $\Delta: X \rar X \times_K X$. Let $\calg{I}$ be
the sheaf of ideals defining the diagonal in $X \times_K X$. \\

\subsection{The complete Bar complex and the complete Hochschild chain complex}

Let $\calg{O}_{\Delta} = \Delta_* \calg{O}_X $. On $X \times_K X$, $\calg{O}_{\Delta}$ has a free $\calg{O}_{X \times_K X}$-module
resolution given by the Bar resolution :\\

$B_n(X) = \calg{O}_X \otimes_K ....... \otimes_K \calg{O}_X (n+2 \text{ times })$.  \\

The $\calg{O}_X \otimes_K \calg{O}_X $ - module structure given by multiplication with the extreme factors. If $U = \text{Spec } R$ is an open
affine subscheme of $X$, then the
differential $d: B_n(R) \rar B_{n-1}(R)$ is given by the formula \\

$d(a_0 \otimes ..... \otimes a_{n+1}) = a_0a_1 \otimes a_2 \otimes ... \otimes a_{n+1} - a_0 \otimes a_1a_2 \otimes ....\otimes a_{n+1} + ......
+ {(-1)}^{n} a_0 \otimes .... \otimes a_na_{n+1} \text{  } a_i \in R$
\\

Let  $\calg{I}_n$  be the kernel of the multiplication map  $B_n(X) \rar \strc$. Let $\hb_n : = \ilim_k B_n/{\calg{I}_n}^k$. Note that the
differential $d$ takes $\calg{I}_n$ into $\calg{I}_{n-1}$. The differential in $B_{\bullet}$ thus extends to a differential in the complete Bar
complex $\hb_{\bullet}$. Note that the complete Bar complex comes equipped with the $\calg{I}$-adic topology.  \\

Yekutieli (Lemma 1.2 of [1]) shows that the complete Bar complex gives us a resolution of $\calg{O}_{\Delta}$ in terms of flat $\calg{O}_{X
\times_K X}$-modules. A resolution of $\Delta^* \calg{O}_{\Delta}$ by flat $\calg{O}_{X}$-modules would thus be given by what is called the
completed Hochschild chain complex of $X$. This complex $\hc_{\bullet}$ is defined by
$$\hc_{\bullet} = \Delta^{-1} \hb_{\bullet} \otimes_{\Delta^{-1} \calg{O}_{X \times_K X}} \calg{O}_{X}$$
The complex $\hc_{\bullet}$ is called the complete Hochschild chain complex of $X$. It has a topology induced by that on $\hb_{\bullet}$. \\

Note that if $U = \text{Spec } R$ is an open affine subscheme  of
$X$, and if $M(R)$ denotes $\Gamma (U, M)$ for any $\mathcal{O}_{X}$
- module $M$, and $I$ is the kernel of the multiplication map
$R^{\otimes_K n+2} \rar R$ we have
$$B_n(R) = R^{\otimes_K n+2} \text{   } \hb_n(R) = \ilim R^{\otimes_K n+2}/I^k$$
$$C_n(R) = R^{\otimes_K n+1}$$ and the differential $d: C_n(R) \rar
C_{n-1}(R)$ is given by
$$d(a_0 \otimes ..... \otimes a_n) = a_0a_1 \otimes ..... \otimes a_n - a_0 \otimes a_1a_2 ..... \otimes a_n + ..... + {(-1)}^n a_n
a_0 \otimes ..... \otimes a_{n-1} $$
$$\hc_n(R) = \ilim R^{\otimes_K n+2}/I^k \otimes_{R^{\otimes_k n+2}} C_n(R)$$ and the differential on $\hc_{\bullet}$ is the one induced
by $d$. \\

\subsection{Hochschild homology, Hochschild cohomology and the completed Hochschild cochain
complex}We begin with the following definitions. \\

{\it Definition 1: } The Hochschild homology of $X$ is defined to be $\rhh(\calg{O}_X
, \Delta^* \calg{O}_{\Delta})$. \\

{\it Definition 2: } The complex of continuous Hochschild cochains on $X$ is the complex \\ $\rhhh_{\calg{O}_X}(\Delta^*
\calg{O}_{\Delta},\calg{O}_X) = \rhhh_{\calg{O}_{X \times_K X}}(\calg{O}_{\Delta},\calg{O}_{\Delta})$. \\

 {\it Fact 1: } Yekutieli (Theorem 0.3 of [1]) shows that
the complex of continuous Hochschild co-chains is given by the complex $\chh_{\calg{O}_{X \times_K X}}(\hb_{\bullet}(X), \calg{O}_{\Delta})$.
This is seen to be equal to the complex $\chh_{\calg{O}_{X}}(\hc_{\bullet}(X),
\calg{O}_{X})$ . Here $\calg{O}_{\Delta}$ and $\strc$ are both given the discrete topology. \\

{\it Fact 2: } It is proven by Yekutieli (Proposition 1.6 of [1]) that the complex of continuous Hochschild cochains on $X$ is none other than
the complex of polydifferential operators on $X$ introduced by Kontsevich [2]. In other words, if $U = \text{Spec} R$ is an open affine
subscheme of $X$, and if ${C_{cd}}^n(X) = \chh(\hc_n(X),\calg{O}_{X})$, then  $${C_{cd}}^n(U) = \{ f \in \text{Hom}_K(R^{\otimes n}, R) | f
\text{ is a differential operator in each factor } \} $$. We shall henceforth denote
${C_{cd}}^n(X)$ by $\dd^n(X)$. \\

{\it Fact 3 :} We recall from Kontsevich (Section 3.4.2) [2]  that
if $U= \text{Spec } R$ and if $f \in {C_{cd}}^n(U)$, then the
differential $d$ of ${C_{cd}}^{\bullet}$ is given by
 $$ df(a_0 \otimes .... \otimes a_{n}) = a_0f(a_1 \otimes ... \otimes a_{n})+ \sum_{i=1}^{i=n} {(-1)}^i f(.. \otimes a_{i-1}a_i \otimes ... \otimes a_n)
+ {(-1)}^{n+1} f(a_0 \otimes...\otimes a_{n-1})a_n $$ for all $a_0,....,a_n \in R$.\\

{\it Remark :} What we refer to here as the complex of polydifferential operators is a shifted version of what Kontsevich [2] refers to as the
complex of polydifferential operators. Kontsevich's complex of polydifferential operators is, in our notation, $\dd^{\bullet}(X)[1]$. \\

\section{Hopf algebra structure on $\dd^{\bullet}(X)$ }

In this section we describe the operations that make $\dd^{\bullet}(X)$ a Hopf algebra in $\bbcc$ . By this, we mean that $\dd^{\bullet}(X)$ has
a multiplication $m$, a comultiplication $\Delta$, a unit $\eta$ and a counit $\epsilon$ all of which are morphisms in $\bbcc$.  \\

Recall that $\dd^1(X)$ is the sheaf of differential operators on $X$. Note that $\dd^1(X)$ is a {\it left} $\strc$ module.
We have the following proposition   \\

\begin{prop}
As $\strc$ modules, $\dd^n(X)$ is isomorphic to $\dd^1(X)^{\tenp n}$.
\end{prop}

\begin{proof}
This is something that can be checked locally. Consider an open affine subscheme $U$ of $X$ with local coordinates $\{x_1,....,x_m\}$ . Then,
an element of $\dd^n(U)$ is given by an map of the form \\
$$f_1 \otimes ... \otimes f_n \leadsto \sum_{(I_1,...,I_n)} C_{(I_1,...,I_n)}(x_1,..,x_m)  \partial_{I_1}f_1.....\partial_{I_n}f_n $$
Here the $I_j$ 's are multi indices and $\partial_{I_j}$ is the partial derivative corresponding to $I_j$ . The above polydifferential operator
maps to $\sum_{(I_1,..,I_n)} C_{(I_1,....,I_n)}(x_1,...,x_m) \partial_{I_1} \otimes ... \otimes \partial_{I_n}$. This gives us a well defined
map from $\dd^n(U)$ to $\dd^1(U)^{\otimes_{{\calg O}_U} n}$. On the other hand, we have a map from $\dd^1(U)^{\otimes_{{\calg O}_U} n}$ to
$\dd^n(U)$ which takes $D_1 \otimes ... \otimes D_n$ to the polydifferential operator $f_1\otimes .... \otimes f_n \leadsto
D_1(f_1)....D_n(f_n)$. These
maps are clearly inverses of each other. \\
\end{proof}

We now describe the Hopf algebra structure on $\dd^{\bullet}(X)$. \\

{\it Multiplication on $\dd^{\bullet}(X)$ : } Let $U = \text{Spec }
R$ be an affine open subscheme of $X$. Let $D_1 \in \dd^k(U)$ and
$D_2 \in \dd^l(U)$. Then we can set $$m(D_1,D_2)(a_1 \otimes ....
\otimes a_{k+l} )= D_1(a_1 \otimes ... \otimes a_k)D_2(a_{k+1}
\otimes... \otimes a_{k+l})$$ for all $a_1,...,a_{k+l} \in R$. This
defines the multiplication $m$ on $\dd^{\bullet}(X)$. Note that \\
$m(D_1,D_2) = D_1 \otimes D_2$
after identifying $\dd^n(X)$ with ${\dd^1(X)}^{\otimes_{\strc} n}$. \\

{\it Comultiplication on $\dd^{\bullet}(X)$ :} On the other hand, if $f_1, ...,f_n$ are differential operators on an open subscheme $U$, then
$f_1......f_n \in \dd^n(U)$ and we can
set \\

$\Delta(f_1....f_n) = \sum_{p + q =n} \sum_{\sigma \text { a (p,q)-shuffle }} \text{sgn}(\sigma) f_{{\sigma}(1)} \otimes ..... \otimes
f_{{\sigma}(p)} \bigotimes f_{{\sigma}(p+1)}
\otimes .... \otimes f_{{\sigma}(p+q)}$. \\

This gives us a well defined map $\Delta: \dd^{\bullet}(U) \rar
\dd^{\bullet}(U) \otimes_{\calg{O}_{X}} \dd^{\bullet}(U)$. This can
be easily seen to commute with restrictions, thus giving us a map \\
$\Delta: \dd^{\bullet}(X) \rar \dd^{\bullet}(X) \otimes
\dd^{\bullet}(X)$, the tensoring being
over $\calg{O}_{X}$. \\

{\it Unit for $\dd^{\bullet}(X)$ :} We have the obvious inclusion map $\eta: \strc \rar \dd^{\bullet}(X)$. On $U= \text{Spec } R$, this is just
the inclusion $R \hookrightarrow \dd^0(U) \hookrightarrow \dd^{\bullet}(U)$. \\

{\it Counit for $\dd^{\bullet}(X)$ :} We also have a projection
$\epsilon: \dd^{\bullet}(X) \rar \strc$. On $U= \text{Spec } R$ with
local coordinates $x_1,..., x_m$ , this takes  \\ $C_0(x_1,...,x_m)
+ \sum_{I_1,...,I_r} C_{I_1,...,I_r}(x_1,...,x_m) \partial_{I_1}
\otimes... \otimes
\partial_{I_r}$ to $C_0(x_1,...,x_m)$. \\

We now have the following fact \\

\begin{prop}

The multiplication $m$, the comultiplication $\Delta$, the unit $\eta$ and the counit $\epsilon$ together make $\dd^{\bullet}(X)$ a Hopf-
algebra in $\bbcc$
\end{prop}

\begin{proof}
Clearly, $\epsilon$ is an algebra homomorphism. Also, $\Delta(id) = id \tenp id $ tells us that $\eta$ is a coalgebra homomorphism. The fact
that $\Delta$ is an algebra homomorphism is exactly analogous to the fact that the comultiplication of the tensor algebra of a vector space over
a field of characteristic $0$ is an algebra homomorphism. Similarly the fact that $m$ is a co-algebra homomorphism is proven in exactly the same way by
which one proves that the product in the graded tensor algebra of a vector space over a field of characteristic $0$ is a co-algebra homomorphism.\\

The only things that remain to be checked are that the differential follows the Leibniz rule and respects co-multiplication - the latter fact
following from the fact that the Hochschild {\it boundary} is a graded derivation with respect to the shuffle product. (see Loday [4]
Proposition 4.2.2). Let $U = \text{Spec R}$ be an open affine subscheme of $X$. Let $a_1,...,a_{k+l+1} \in R$. Let $D_1 \in \dd^k(U)$ and $D_2
\in \dd^l(U)$. The following calculation verifies that the differential $d$ on $\dd^{\bullet}(X)$ obeys the Leibniz rule with respect to
the multiplication $m$. \\

$$d(D_1 \otimes D_2)(a_1 \otimes .... \otimes a_{k+l+1}) = a_1 (D_1 \otimes D_2)(a_2 \otimes .... \otimes a_{k+l+1}) - (D_1 \otimes D_2)(a_1a_2
\otimes ... \otimes a_{k+l+1}) + ... $$ $$... + {(-1)}^{k}(D_1
\otimes D_2)(a_1 \otimes ... \otimes a_ka_{k+1} \otimes ... \otimes
a_{k+l+1}) +...+{(-1)}^{k+l+1} (D_1 \otimes D_2) (a_1 \otimes..
\otimes a_{k+l}) a_{k+l+1} $$

$$ = a_1D_1(a_2 \otimes ... \otimes a_{k+1}).D_2(a_{k+2} \otimes ... \otimes a_{k+l+1})  - D_1(a_1.a_2 \otimes ... \otimes a_{k+1})D_2(a_{k+2} \otimes ... \otimes a_{k+l+1})
+.. $$ $$. +{(-1)}^k D_1(a_1 \otimes .. \otimes a_ka_{k+1})
D_2(a_{k+2} \otimes .. \otimes a_{k+l+1}) $$ $$+ {(-1)}^{k+1}
D_1(a_1 \otimes .. \otimes a_k)a_{k+1}D_2(a_{k+2} \otimes .. \otimes
a_{k+l+1}) $$

$$ + {(-1)}^k D_1(a_1 \otimes ... \otimes
a_k)a_{k+1}D_2(a_{k+2} \otimes .... \otimes a_{k+l+1}) $$ $$ +
{(-1)}^{k+1}D_1(a_1 \otimes ... \otimes a_k)D_2(a_{k+1}a_{k+2}
\otimes ... \otimes a_{k+l+1}) +  $$ $$.... +{(-1)}^{k+l+1}D_1(a_1
\otimes ... \otimes a_k)D_2(a_{k+1} \otimes ... \otimes a_{k+l})
a_{k+l+1} $$

$$ = dD_1(a_1 \otimes ... \otimes a_{k+1})D_2(a_{k+2} \otimes ... \otimes a_{k+l+1}) + {(-1)}^k D_1(a_1 \otimes ... \otimes a_{k})dD_2(a_{k+1} \otimes ... \otimes
a_{k+l+1})$$

\end{proof}

\begin{cor}
The maps in $\dcat$ induced by $m$, $\Delta$ $\eta$ and $\epsilon$ make $\dd^{\bullet}(X)$ a Hopf-algebra in $\dcat$.

\end{cor}

Remark: In fact, the Hopf algebra structure on $\dd^{\bullet}(X)$ is that of the graded tensor algebra. As $\strc$ - modules $\dd^n(X)$ is
isomorphic to ${\dd^1(X)}^{\tenp n}$ (Proposition 1) .  Thus, as far as the Hopf algebra structure is concerned, $\dd^{\bullet}(X)$ is
isomorphic to the tensor algebra $T(\dd^1(X))$ generated over $\strc$ by $\dd^1(X)$ in degree $1$. They are isomorphic as Hopf algebras in $Gr
(\text{mod}-\strc)$. But there is a nontrivial differential (Hochschild cochain differential) on
$\dd^{\bullet}(X)$. \\

\section{A decomposition of $\dd^{\bullet}(X)$}

\subsection{The decomposition}

Recall that $L(\dd^1(X))$ denotes the free Lie algebra generated over $\strc$ by $\dd^1(X)$ concentrated in degree $1$. As graded $\strc$
modules, it is a submodule of $\dd^{\bullet}(X)$. Moreover,

\begin{prop}

The differential on $\dd^{\bullet}(X)$ preserves $L(\dd^1(X))$ \\

\end{prop}

\begin{proof}
By proposition 2, the differential $d$ on $\dd^{\bullet}(X)$ obeys Leibniz rule with respect to the multiplication $m$ on $\dd^{\bullet}(X)$. It
follows that it obeys Leibniz rule with respect to the Lie bracket $[,]$ induced by $m$ on $\dd^{\bullet}(X)$. The restriction of $[,]$ to
$L(\dd^1(X))$ is precisely the Lie bracket on $L(\dd^1(X))$ and will also be denoted by $[,]$. Since $\dd^1(X)$ (in degree $1$) generates
$L(\dd^1(X))$ as a Lie algebra over $\strc$, and $d$ obeys the Leibniz rule with respect to $[,]$, it is enough to check that $d(\dd^1(X))$ is
contained in $L(\dd^1(X))$. \\

Since the differential $d$ and the bracket $[,]$ are $\strc$-linear, we only need to check that if $U = \text{Spec } R$ with local coordinates
$x_1,...,x_m$ and if $I$ is a multi-index,and if $\partial_I$ denotes the corresponding partial derivative, then
$$d(\partial_I(a_1 \otimes a_2)) \in L(\dd^1(U))$$ Recall that if $I= (p_1,...p_m)$ and if $J= (q_1,...,q_m)$ then $ J \prec I$ if $q_i \leq p_i \text{ } \forall i$
 and $J \neq I$. Let $I-J:= (p_1-q_1,...,p_m-q_m)$.Now, $$d(\partial_I(a_1 \otimes a_2)) = a_1 \partial_I a_2 -\partial_I(a_1.a_2) + a_2.
\partial_I(a_1) $$ $$= - \frac{1}{2} \sum_{J \prec I} C_{IJ} (\partial_J(a_1). \partial_{I-J}(a_2) + \partial_{I-J}(a_1)\partial_J(a_2)$$
$$= - \frac{1}{2} \sum_{J \prec I} C_{IJ} [\partial_J, \partial_{I-J}](a_1 \otimes a_2)$$. Here $C_{IJ}$ are some rational constants.
This completes the desired verification.\\

\end{proof}

{\it The symmetrization map $I$ : }We have a symmetrization map $$I:\oplus_k \sss^k(L(\dd^1(X))) \rar \dd^{\bullet}(X)$$ Let $U = \text{Spec }
R$ be an open affine subscheme of $X$, and let $z_1,...,z_k$ be homogenous elements of $L(\dd^1(U))$ of degrees $d_1,...,d_k$ respectively. If
$s(\sigma)$ is the sign such that $z_1....z_k =s(\sigma) z_{\sigma(1)}.... z_{\sigma(k)} \in \sss^k(L(\dd^1(U))) $ then
$$I(z_1....z_k) = \frac{1}{k!} \sum_{\sigma \in S_k} s(\sigma) z_{\sigma(1)} \otimes ... \otimes z_{\sigma(k)} \in \dd^{\bullet}(X) $$ Note that $s(\sigma)$
depends on $d_1,...,d_k$ and $\sigma$. \\

We now have the following key lemma \\

\begin{lem}

The symmetrization map $I$ is an isomorphism in $\bbcc$. \\

\end{lem}

\begin{proof}
It follows directly from Propositions 2 and 3 and the definition of $I$ that $I$ is a map of complexes of $\strc$ modules. It therefore suffices
to show that $I$ is a map of graded $\strc$-modules. This can be verified locally. \\

Note that $\dd^1(X)$ is locally free. Let $U = \text{Spec } R$ be an affine open subscheme of $X$ such that $\dd^1(U)$ is trivial on $U$. Then,
$\dd^1(U) = V \otimes_K \calg{O}_{U}$. If $L(V)$ is the free Lie algebra generated over $K$ by $V$ in degree $1$ and if $T(V)$ is the tensor
algebra of $V$, then  $L(\dd^1(U)) = L(V) \otimes_K \calg{O}_{U}$ and $\dd^{\bullet}(U) = T(V) \otimes_K \calg{O}_{U}$ as graded
$\calg{O}_{U}$-modules. \\

Let $I_V:\oplus_k \sss^k(L(V)) \rar T(V)$ be the symmetrization map. Let  $L_1,...,L_k$ be homogenous elements of $L(V)$ of degrees
$d_1,...,d_k$ respectively. Let $s(\sigma)$ be the sign such that $L_1...L_k = s(\sigma) L_{\sigma(1)}...L_{\sigma(k)}$ in $\sss^k(L(V))$. Then,
 \\ $I_V (L_1...L_k)= \frac{1}{k!} \sum_{\sigma \in S_k} s(\sigma) L_{\sigma(1)} \otimes... \otimes L_{\sigma(k)}$. Note that \\ $\oplus
\sss^k(L(\dd^1(U))) = \oplus_k \sss^k(L(V)) \otimes_K \calg{O}_{U}$ and that $I = I_V \otimes id_{\calg{O}_{U}}$. It thus suffices to show that
$I_V$ is an
isomorphism of graded $K$-vector spaces. This is Proposition 17 of the Appendix to this paper. \\
\end{proof}

\subsection{The Hochcshild-Kostant-Rosenberg map}

Throughout this subsection let \\ $U = \text{Spec R}$ be an affine
open subscheme of $X$ with local coordinates $x_1,...,x_m$. Recall
(Yekutieli [1] Theorem 4.8) that the Hochschild-Kostant-Rosenberg
map \\ $\hkr: \oplus_i \texal{X}{i} \rar \dd^{\bullet}(X)$ is a map
of complexes which is a quasi isomorphism. Also recall (Yekutieli
[1] page 14) that if  $a_1...,a_k \in R$ then
$$\hkr(\frac{\partial}{\partial x_{i_1}} \wedge ... \wedge
\frac{\partial}{\partial x_{i_k}})(a_1 \otimes .... \otimes a_k) =
\frac{1}{k!} \sum_{\sigma \in S_k} \text{sgn}(\sigma) \frac{\partial
a_1}{\partial x_{i_{\sigma(1)}}}......\frac{\partial a_k}{\partial
x_{i_{\sigma(k)}}} $$

Note that we have a map of complexes $\beta: T_X[-1] \rar
L(\dd^1(X))$ such that on $U$,  $\beta(\frac{\partial}{\partial
x_i})(a) = \frac{\partial a}{\partial x_i}  \text{ } \forall a \in R
$. Note that $\beta$ induces a map of complexes \\  $\sss^{\bullet}
\beta : \oplus_i
\texal{X}{i} \rar \oplus_k \sss^k(L(\dd^1(X)))$. \\

\begin{prop} $\sss^{\bullet} \beta$ is a quasi-isomorphism. \end{prop}

\begin{proof}
Observe that if $a_1,...,a_k \in R$ then $$I \circ \sss^{\bullet} \beta (\frac{\partial}{\partial x_{i_1}} \wedge ... \wedge
\frac{\partial}{\partial x_{i_k}})(a_1 \otimes .... \otimes a_k) = \frac{1}{k!} \sum_{\sigma \in S_k} \text{sgn}(\sigma) \frac{\partial
a_1}{\partial x_{i_{\sigma(1)}}}......\frac{\partial a_k}{\partial x_{i_{\sigma(k)}}}
$$. Thus $I \circ \sss^{\bullet} \beta  = \hkr$.   \\

Now $\hkr$ is a quasi-isomorphism (by Yekutieli [1] theorem 4.8) and  $I$ is an
isomorphism of complexes of $\strc$ modules by Lemma 1 and therefore a quasi-isomorphism. Thus, $\sss^{\bullet} \beta$ is a quasi-isomorphism. \\
\end{proof}

\begin{prop} $\sss^k \beta : \texal{X}{k} \rar \sss^k(L(\dd^1(X)))$ is a quasi-isomorphism. \end{prop}

\begin{proof} Given a complex $\calg M$ in $\bbcc$, let ${\calg H}^*(\calg M)$ denote its cohomology , which is a graded $\strc$ module.
Then ${\calg H}^* (\texal{X}{k}) = \texal{X}{k}$ as $\texal{X}{k}$
has zero differential. Note that $\sss^{\bullet} \beta$ induces a
map ${\sss^{\bullet} \beta}_*: \oplus_k \texal{X}{k} \rar {\calg
H}(\oplus_k \sss^k(L(\dd^1(X)))) $ which is an isomorphism of graded
$\strc$-modules by Proposition 4. But $\sss^{\bullet} \beta =
\oplus_k \sss^k \beta$ by definition. Thus, ${\sss^{\bullet}
\beta}_* = \oplus_k {\sss^k \beta}_* $ where  \\ ${\sss^k \beta}_*:
\texal{X}{k} \rar {\calg H} (\sss^k(L(\dd^1(X))))$ is the map
induced on cohomology by $\sss^k \beta$. Since ${\sss^{\bullet}
\beta}_*$ is an isomorphism of graded $\strc$-modules, it follows
that for all $k$, ${\sss^k \beta}_*$ is an
isomorphism of graded $\strc$-modules.\\
\end{proof}

In particular $\beta: T_X[-1] \rar L(\dd^1(X))$ is a quasi-isomorphism. We state this as a separate corollary in order to highlight it. Thus,\\

\begin{cor} $\beta :T_X[-1] \rar L(\dd^1(X))$ is a quasi-isomorphism \end{cor}

\section{Theorem 1}

\subsection{Precise statement and proof of Theorem 1}

Throughout this section let \\ $U = \text{Spec R}$ be an arbitrary
affine subscheme of $X$. Let $z_1,...,z_k,y$ be homogenous elements
of $L(\dd^1(U))$ of degrees $d_1,..,d_k,d$ respectively.  Let
$$\omega: \oplus_k \sss^k(L(\dd^1(X))) \otimes L(\dd^1(X)) \rar \oplus_k \sss^k(L(\dd^1(X))) \otimes L(\dd^1(X))$$ be the morphism in $\bbcc$
such that on $U$, $$\omega(z_1...z_k \otimes y) = \sum_{i=1}^{i=k}
{(-1)}^{d_i(d_{i+1}+...+d_k)} z_1..\widehat{z_i} ... z_k \otimes
[z_i, y]$$ Note that if $\mu$ denotes the multiplication on
$\sss^{\bullet}(L(\dd^1(X)))$, then $\mu \circ \omega = \adj$ where
$\adj$ denotes the right adjoint
action of $L(\dd^1(X))$ on $\sss^{\bullet}(L(\dd^1(X)))$. Let $I$ be as in Section 4. \\

More generally , if $M$ is a locally free $\strc$ module, let $T(M)$
denote the (graded) tensor algebra of $M$. Let $L(M)$ denote the
free Lie algebra generated over $\strc$ by $M$ in degree $1$. On an
open affine subscheme  $U = {\text{Spec }}R$ , let $q_1,...,q_k,r$
be homogenous elements of $L(M|_U)$ of degrees $d_1,...,d_k, d$
respectively. Let $$\omega_M: \sss^{\bullet}(L(M)) \otimes L(M) \rar
\sss^{\bullet}(L(M)) \otimes L(M) $$ be the morphism of $\strc$
modules such that $$\omega_M(q_1...q_k \otimes r) = \sum_{i=1}^{i=k}
{(-1)}^{d_i(d_{i+1}+...+d_k)} q_1...\widehat{q_i}...q_k \otimes
[q_i,r] $$ Let $I_M: \sss^{\bullet}(L(M)) \rar T(M)$ denote the
symmetrization map such that $$I_M(q_1...q_k) = \sum_{\sigma \in
S_k} s(\sigma) q_{\sigma(1)} \otimes .... \otimes q_{\sigma(k)} $$
where $s(\sigma)$ is the sign such that $q_1...q_k =
s(\sigma) q_{\sigma(1)}...q_{\sigma(k)}$ in $\sss^k(L(M))$. \\

If $V$ is a vector space over $K$, let $T(V)$ be the (graded) tensor algebra generated by $V$ over $K$. Let $L(V)$ denote the free Lie algebra
generated over $K$ by $V$ in degree $1$. Let $v_1,...,v_k,w$ be homogenous elements of $L(V)$ of degrees $d_1,...,d_k,d$ respectively. Let
$$ \omega_V: \sss^{\bullet}(L(V)) \otimes L(V) \rar \sss^{\bullet}(L(V)) \otimes L(V) $$ be the map such that $$\omega_V(v_1...v_k \otimes w) =
 \sum_{i=1}^{i=k} {(-1)}^{d_i(d_{i+1}+...+d_k)}
v_1...\widehat{v_i}...v_k \otimes [v_i,w] $$ Let $I_V:
\sss^{\bullet}(L(V)) \rar T(V)$ denote the symmetrization map such
that $$I_V(v_1...v_k) = \sum_{\sigma \in S_k} s(\sigma)
v_{\sigma(1)} \otimes .... \otimes v_{\sigma(k)} $$ where
$s(\sigma)$ is the sign such that $v_1...v_k =
s(\sigma) v_{\sigma(1)}...v_{\sigma(k)}$ in $\sss^k(L(V))$. \\

\begin{thm1}

The following diagram commutes in $\bbcc$ \\

$$\begin{CD}
\dd^{\bullet}(X) \otimes L(\dd^1(X)) @> m >> \dd^{\bullet}(X) \\
@AA{I \otimes id}A           @A{I}AA \\
\sss^{\bullet}(L(\dd^1(X))) \otimes L(\dd^1(X)) @> \mu \circ \frac{\omega}{1-e^{-\omega}}>> \sss^{\bullet}(L(\dd^1(X)))\\
\end{CD} $$

\end{thm1}

\begin{proof}

{\it Step 1 } (Checking that all morphisms involved commute with the relevant differentials) \\

Before we proceed, we note that the product $m$ on $\dd^{\bullet}(X)$ and $I$ commute with the relevant differentials by Proposition 2 and
Lemma 1 respectively. \\

Further , we need to see that $\mu \circ \frac{\omega}{1-e^{-\omega}}$ commutes with the relevant differentials. Since $\mu$ commutes with the
relevant differentials, we only need to check that $\frac{\omega}{1-e^{-\omega}}$ commutes with the relevant differentials. The latter
expression is a power series in $\omega$. It is therefore enough to verify that $\omega$ commutes with the relevant differentials. This only
needs to be checked locally. The following calculations are done to complete the check \\

$$\omega(z_1z_2.....z_k \otimes y) = \sum_i {(-1)^{d_i(d_{i+1}+....+ d_k)}} z_1.....\widehat{z_i} ...z_k \otimes [z_i,y] $$
$$d(\omega(z_1z_2.....z_k \otimes y)) = \sum_i {(-1)^{d_i(d_{i+1}+....+d_k)}} ( \sum_{j \neq i} {(-1)^{d_1+..\widehat{d_i} ...+d_{j-1}}}
z_1...dz_j....\widehat{ z_i} ...z_n \otimes [z_i,y]$$ $$ +
{(-1)^{d_1+... \widehat{d_i} ..+d_k}} z_1..... \widehat{z_i} ...z_k
\otimes d([z_i,y]) $$
$$ d( z_1z_2.....z_k \otimes y) = (\sum_j {(-1)^{d_1+...+d_{j-1}}} z_1...dz_j....z_k \otimes y) + {(-1)^{\sum d_i}} z_1......z_k \otimes dy $$
$$\omega(d(z_1z_2.....z_k \otimes y)) = \sum_j\{ ( \sum_{i < j}  {(-1)^{d_1+...+d_{j-1}}}  {(-1)^{d_i(d_{i+1}+...+d_{k}+1)}} z_1...\widehat{z_i}..
 dz_j.... ...z_k \otimes [z_i,y]
 $$ $$+ \sum_{i>j} {(-1)^{d_1+...+d_{j-1}}}  {(-1)^{d_i(d_{i+1}+...+d_{k})}} z_1...dz_j....\widehat{z_i }...z_k \otimes [z_i,y]) $$ $$ +
 {(-1)^{d_1+...+d_{j-1}}} {(-1)^{(d_j+1)(d_{j+1}...d_k)}} z_1..... \widehat{z_j} ...z_k \otimes [dz_j,y] $$ $$ +
{(-1)^{\sum d_i}} {(-1)^{d_j(d_{j+1}+...+d_k)}}  z_1.....
\widehat{z_j} ...z_k \otimes [z_j, dy]\} $$

We now compare the coefficients of $[z_i, y]$, $[dz_i,y]$ and $[z_i,dy]$ in $d(\omega(z_1....z_k \otimes y))$ and $\omega(d(z_1....z_k \otimes
y))$
and check that they are equal . \\

{\it Step 2} (Reduction to an analogous result for the graded free Lie algebra generated by a vector space over a field of characteristic $0$)\\

Having checked that all morphisms in the diagram given in the proposition are morphisms in $\bbcc$ , it suffices to verify that the above
diagram commutes in the category of graded $\strc$-modules. In other words, we can "forget the differentials involved". Recall that $\dd^1(X)$
is a locally free $\strc$ module. It
therefore suffices to prove that if $M$ is a  locally free $\strc$ module, the following diagram commutes in the category of graded $\strc$ modules. The map $m_M
$ in the diagram below is the multiplication in $T(M)$.\\

$$\begin{CD}
T(M) \otimes L(M) @> m_M >> T(M) \\
@AA{I_M \otimes id}A @A{I_M}AA \\
\sss^{\bullet}(L(M)) \otimes L(M) @> \mu \circ \frac{\omega_M}{1-e^{-\omega_M}}>> \sss^{\bullet}(L(M)) \\
\end{CD} $$

All morphisms in the diagram in this proposition are $\strc$ module homomorphisms. It therefore suffices to check the claim that the above
diagram commutes locally. We may therefore, without loss of generality assume that $M$ is a free $\strc$ module i.e, $M = \strc \otimes_K V$ for
some $K$ vector space $V$. Then, $T(M) = \strc \otimes_K T(V)$ , $L(M) = \strc \otimes_K L(V)$ and $\sss^{\bullet}(L(M)) = \strc \otimes_K
\sss^{\bullet}(L(V))$.  \\

Since the morphisms in the commutative diagram before the previous paragraph are all $\strc$ linear, it suffices to check that the following
diagram commutes in the category of graded $K$ vector spaces. The map $m_V$ in the diagram below is the multiplication in $T(V)$. \\

$$\begin{CD}
T(V) \otimes L(V) @> m_V >> T(V) \\
@AA{I_V \otimes id}A @A{I_V}AA \\
\sss^{\bullet}(L(V)) \otimes L(V) @> \mu \circ \frac{\omega_V}{1-e^{-\omega_V}}>> \sss^{\bullet}(L(V)) \\
\end{CD} $$

This is Theorem 6 of the Appendix. \\

\end{proof}

\subsection{Some remarks on Theorem 1}

Let $A$ be an associative algebra in $\bbcc$. In other words, there is a multiplication morphism $\mu_A: A \otimes A \rar A$ in $\bbcc$ which is
associative i.e, $\mu_A \circ (\mu_A \otimes id)= \mu_A \circ (id \otimes \mu_A)$ as morphisms in $\bbcc$ from $A \otimes A \otimes A$ to $A$.\\
\\
A Lie algebra $L$ in $\bbcc$ is an object in $\bbcc$ equipped with a morphism $[,]_L: L \otimes L \rar L$ in $\bbcc$ such that \\
(i) $[,]_L = - [,]_L \circ \tau$ where $\tau: L \otimes L \rar L \otimes L$ is the swap map. \\
(ii) $[,]_L \circ (id \otimes [,]_L) =[,]_L \circ ([,]_L \otimes id) + [,]_L \circ (id \otimes [,]_L) \circ (\tau \otimes id) $ \\

Note that any algebra $A$ in $\bbcc$ has a Lie algebra structure with  \\ $[,]_A = \mu_A \circ (id -\tau)$. \\

Note that by Proposition 2 , $\dd^{\bullet}(X)$ is an algebra in $\bbcc$ and $L(\dd^1(X))$ is a Lie algebra in $\bbcc$. \\

Given a Lie algebra $L$ in $\bbcc$, its universal enveloping algebra
(if it exists) is an algebra $U(L)$ in $\bbcc$ together with a
morphism $i: L \rar U(L)$ of Lie algebras such that given any
morphism $f: L \rar A$ of Lie algebras from $L$ to an algebra $A$ in
$\bbcc$, there exists a unique morphism $\bar{f}:U(L) \rar A$ of
algebras in $\bbcc$ such that $f = \bar{f} \circ i$.  We now prove
the
following lemma \\

 \begin{lem}

 $\dd^{\bullet}(X)$ is the universal enveloping algebra of $L(\dd^1(X))$ in $\bbcc$.

 \end{lem}

 \begin{proof}

Let $A$ be an algebra in $\bbcc$. Let $f: L(\dd^1(X)) \rar A$ be a morphism of Lie algebras in $\bbcc$. In other words, $f \circ [,]_L = \mu_A
\circ (id -\tau) \circ (f \otimes f)$ as morphisms in $\bbcc$. Let $A^n$ denote the degree $n$ term of the complex $A$. Let $f^n: L(\dd^1(X))
\cap \dd^n(X) \rar A^n$ be the degree $n$ component of the morphism $f$. In particular, $f^1: \dd^1(X) \rar A^1$. \\

Note that ${f^1}^{\otimes n} :\dd^n(X) = \dd^1(X)^{\otimes n} \rar {(A^1)}^{\otimes n}$ is a map of $\strc$-modules. Note that the $n$-fold
multiplication $\mu_{n,A}:= \mu_A \circ (\mu_A \otimes id) \circ .... \circ (\mu_A \otimes id^{\otimes n-1}):A^{\otimes n} \rar A$ maps
${(A^1)}^{\otimes n}$ to $A^n$. Set $\bar{f}^n: \dd^n(X) \rar A^n$ to be the composite $\mu_{n,A} \circ {f^1}^{\otimes n}$. \\

Let $\bar{f}: \dd^{\bullet}(X) \rar A$ be the map of graded $\strc$-modules whose degree $n$ component is $\bar{f}^n$. We need to check that
$\bar{f}$ is indeed a map of complexes of $\strc$-modules. This can be checked locally. Suppose $U = \text{Spec } R$ is an affine open subscheme
of $X$, and if $D_1,...,D_n$ are differential operators on $U$. If $d$ and $d_A$ denote the differentials on $\dd^{\bullet}(U)$ and $A |_U$
respectively, then $$d_A(f^1(D_1)...f^1(D_n)) = \sum_{i=1}^{i=n} {(-1)}^{i-1} f^1(D_1)...d_A(f^1(D_i))...f^1(D_n) $$ $$ = \sum_{i=1}^{i=n}
{(-1)}^{i-1} f^1(D_1)...f^2(d(D_i))...f^1(D_n) $$ The last equality is because $f$ is a map of complexes. \\

Let $ L(\dd^1(U))^k$ denote the degree $k$ term of the complex $L(\dd^1(U))$. Note that $d(D_i) \in L(\dd^1(U))^2$ by Proposition 3. Also, if
$[,]$ denotes the bracket on $L(\dd^1(U))$ , then $f^2 \circ [,] = \mu_A \circ (id- \tau) \circ f^1 \otimes f^1 $ since $f$ is a Lie algebra
homomorphism in $\dcat$. Also, $(id -\tau) \circ (f^1 \otimes f^1) = (f^1 \otimes f^1) \circ (id - \tau)$.   Moreover, $[,]: \dd^1(U)^{\otimes
2} \rar L(\dd^1(U))^2 $ is surjective. It follows that $f^2(d(D_i)) = (f^1 \otimes f^1) d(D_i)$ ($d(D_i)$ on the right hand
side is thought of as an element of $\dd^2(U)$ after identifying $L(\dd^1(U))^2$ with its image in $\dd^2(U)$ ). \\

It follows that $$d_A(f^1(D_1)...f^1(D_n)) =  \sum_{i=1}^{i=n} {(-1)}^{i-1} f^1(D_1)...f^2(d(D_i))...f^1(D_n)$$ $$ =   \sum_{i=1}^{i=n}
{(-1)}^{i-1} f^1(D_1)...(f^1 \otimes f^1) d(D_i))...f^1(D_n) $$ $$ = \sum_{i=1}^{i=n} {f^1}^{\otimes n+1} (D_1 \otimes ... \otimes d(D_i)
\otimes ... \otimes D_n) $$ $$ = \bar{f}^{n+1} d(D_1 \otimes ... \otimes D_n)$$ This shows that $\bar{f}$ is indeed a map of complexes. That
$\bar{f}$ is a map of algebras in $\bbcc$ is now immediate from its construction. \\

If $i: L(\dd^1(X)) \rar \dd^{\bullet}(X)$ is the restriction of $I$ to $L(\dd^1(X))$ we need to check that $\bar{f} \circ i = f$. Since the maps
involved are maps of complexes of $\strc$-modules, it is enough to check that $\bar{f} \circ i = f$ as maps of graded $\strc$-modules. This can
again be checked locally. \\

Let $\tau_k$ denote the $k$-cycle $(n-k+1 \text{ } n-k+2 .... \text{ } n)$ of $S_n$. Consider the element $\sigma_n:= (1 - \tau_n) ...
(1-\tau_2)$ of group ring $KS_n$ of $S_n$. Recall that $S_n$ acts on ${A^1}^{\otimes n}$ on the right by $a_1 \otimes ... \otimes a_n \leadsto
\text{sgn}(\sigma) a_{\sigma(1)} \otimes ... \otimes a_{\sigma(n)}$. By $K$-linearity, this extends to an action of $KS_n$ on ${A^1}^{\otimes n}$. \\

Suppose $U = \text{Spec } R$ is an affine open subscheme of $X$, and if $D_1,...,D_n$ are differential operators on $U$. Let $[,]_n:\dd^1(U)
\rar L(\dd^1(U))$ denote the "$n$-fold bracket" i.e, the map taking $D_1 \otimes ... \otimes D_n$ to $[D_1,[D_2[...[D_{n-1},D_n]]]]$. The
following
diagram commutes in the category of graded $R$-modules since $f$ is a Lie algebra homomorphism in $\bbcc$. \\

$$\begin{CD}
\dd^1(U)^{\otimes n} @> {{f^1}^{\otimes n}} >> {A^1}^{\otimes n} \\
@V{[,]_n}VV       @VV{\mu_{n,A} \circ \sigma_n}V \\
L(\dd^1(U))^n @>{f^n}>> A^n \\
\end{CD} $$

Note that $S_n$ acts on the right on $\dd^1(U)^{\otimes n}$ as well.
Also, ${f^1}^{\otimes n} \circ \sigma_n = \sigma_n \circ
{f^1}^{\otimes n}$. Also $i \circ [,]_n = \sigma_n:
\dd^1(U)^{\otimes n} \rar \dd^1(U)^{\otimes n} = \dd^n(U)$. It
follows that \\ $f^n = {f^1}^{\otimes n} \circ i : L(\dd^1(U))^n
\rar A^n$. This proves that $\bar{f} \circ i = f$. \\

Finally, we need to prove that $\bar{f}$ is the unique map with the
required properties. Suppose that if $g:\dd^{\bullet}(X) \rar A$ is
a morphism of algebras in $\bbcc$ such that $g \circ i = f$. Then ,
the restriction of $g$ to $\dd^1(X)$ is precisely $f^1$. It then
follows from the fact that $g$ is an algebra morphism in $\bbcc$
that the restriction of $g$ to $\dd^n(X) = {\dd^1(X)}^{\otimes n}$
is precisely $\bar{f}^n$. This proves that $ g = \bar{f}$. \\

\end{proof}

\subsubsection{Meaning of Theorem 1}

Let $\mathfrak{g}$ denote a finite dimensional Lie algebra over a
field $K$ of characteristic $0$. Let $U(\mathfrak{g})$ denote the
universal enveloping algebra of $\mathfrak{g}$. Let $\pbw:
\sss^{\bullet} \mathfrak{g} \rar U( \mathfrak{g})$ be the
symmetrization map. If $g_1,.. g_k \in \mathfrak{g}$ then,
\\ $\pbw(g_1...g_k) = \sum_{\sigma \in S_k} \frac{1}{k!}
g_{\sigma(1)}*....*g_{\sigma(k)}$ where $*$ denotes
the multiplication in $U( \mathfrak{g})$. Let $\exx: \mathfrak{g} \rar U( \mathfrak{g})$ denote the exponential map. \\

Consider the following calculation of $d(\exx)^{-1}$. Note first
that \\ $\exx(v) = \pbw \circ e^v \text{ } \forall \text{ } v \in
\mathfrak{g}$, where $e^v = 1 + v + \frac{v^2}{2!} + .... \in
\sss^{\bullet} \mathfrak{g}$. Also, \\ $d(\exx)^{-1}_{\exx(X)}(Y) =
\frac{d}{dt} |_{t=0} \exx^{-1}(\exx(X)*\exx(tY))$ . \\ Thus
$\exx^{-1}(\exx(X)*\exx(ty)) = X + d(\exx)^{-1}_{\exx(X))}(Y).t +
h.o.t $. It follows that \\ $\pbw^{-1}(\pbw (e^X) * \pbw(e^{tY})) =
e^{X + d(\exx)^{-1}_{\exx(X)}(Y).t + h.o.t}$. Taking the derivative
with respect to $t$ at $t=0$ on both sides of the previous equation,
we get \\ $\pbw^{-1}(\pbw(e^X) * Y) =
e^X.d(\exx)^{-1}_{\exx(X))}(Y)$ where the multiplication on the
right is the product in $\sss^{\bullet} {\mathfrak{g}}$. \\

Another way of looking at this phenomenon is to say that the calculation of $d(\exx)^{-1}$ is equivalent to specifying
$\varphi$ in the following commutative diagram \\

$$ \begin{CD} U{\mathfrak{g}} \otimes \mathfrak{g} @> \mu >>
U{\mathfrak{g}} \\
@A \pbw \otimes id AA              @AA \pbw A \\
\sss^{\bullet}{\mathfrak{g}} \otimes \mathfrak{g} @> \varphi >>
\sss^{\bullet}{\mathfrak{g}} \\
\end{CD} $$ \\

 Let $$\omega_{
\mathfrak{g}}:\sss^{\bullet}{\mathfrak{g}} \otimes {\mathfrak{g}} \rar \sss^{\bullet}{\mathfrak{g}} \otimes {\mathfrak{g}}$$ be the map such
that for $g_1,...,g_k, y  \in {\mathfrak{g}}$, $$\omega_{\mathfrak{g}}(g_1...g_k \otimes y ) = \sum_i g_1...\widehat{i}..g_k \otimes [g_i,y]$$ .
Note that if $\mu$ is the natural product in $\sss^{\bullet}({\mathfrak{g}})$ then $\mu \circ \omega_g = \adj$ where $\adj$ denotes the right
adjoint action of ${\mathfrak{g}}$ on $\sss^{\bullet} {\mathfrak{g}}$. In the classical situation we know that $d(\exx)^{-1}_{\exx(X))}(Y) =
\frac{\adj(X)}{1-e^{-\adj(X)}} Y$ .From this it follows that $\varphi = \mu \circ \frac{\omega_{\mathfrak{g}}}{1-e^{-\omega_{\mathfrak{g}}}}$. \\

This in short describes how $\pbw$ fails to commute with
multiplication in the classical situation. By Lemma 2,
$\dd^{\bullet}(X)$ is the universal enveloping algebra of
$L(\dd^1(X))$ in $\bbcc$. $I$ is the symmetrization map, and is the
direct analog in our situation of $\pbw$. Theorem 1 , therefore,
specifies how $I$ fails to commute with multiplication and says that
the error term measuring this failure has the same
"$d(\text{exp}^{-1})$ like" form as the corresponding error term in
the classical situation. Since measuring how $\pbw$ fails to commute
with multiplication calculates
$d(\exx)^{-1}$ in the classical situation, we can call Theorem 1 the calculation of $d(\exx)^{-1}$ for the Lie algebra $L(\dd^1(X))$ of $\bbcc$. \\

\section{The Atiyah class of $T_X$}

Recall that $\dcat$ denotes the derived category of bounded below
complexes of $\strc$-modules with coherent cohomology. Let $E$ be a
vector bundle on $X$, and let $J_1(E)$ be the bundle of first jets
of $E$.Recall that the Atiyah class of $E$ is the element in
$\text{Hom}_{\dcat}( E, E \tenp \Omega [1])$ arising out of the
exact sequence $0 \rar E \tenp \Omega \rar J_1(E) \rar E \rar 0$. We
denote the
Atiyah class of $E$ by $\alpha(E)$. \\

We can extend the definition of the Atiyah class of a vector bundle to elements of $\dcat$ given by classes of perfect complexes of vector
bundles on $X$. If $E^{\bullet}$ is such a complex, we have an exact sequence of complexes of vector bundles $$\begin{CD} 0 @>>> E^{\bullet}
\tenp \Omega @>f>> J_1(E^{\bullet}) @>g>> E^{\bullet} \rar 0 \end{CD}$$ We thus have a quasi-isomorphism of complexes $$ Q: \text{tot}(0 \rar
E^{\bullet} \tenp \Omega \rar J_1(E^{\bullet})) \rar E^{\bullet}$$ such that the map  $J_1(E^n) \oplus E^{n-1} \tenp \Omega
 \rar E^n $ is the composite $$\begin{CD} J_1(E^n) \oplus E^{n-1} @>>> J_1(E^n) @>{g^n}>> E^n \end{CD}$$ where the first arrow is projection to the first
factor.  On the other hand, we have a map of complexes
 $$ R: J_1(E^{\bullet}) \oplus E^{\bullet} \tenp \Omega [1]  \rar  E^{\bullet} \tenp \Omega [1]
 $$ given by projection to the second factor.  Consider the element $$R \circ Q^{-1} \in \text{Hom}_{\dcat}(E^{\bullet}, E^{\bullet}
 \tenp \Omega [1])$$ One checks that replacing $E^{\bullet}$ by a complex
 quasi-isomorphic to $E^{\bullet}$ does not give us a different element in
 $\text{Hom}_{\dcat}(E^{\bullet}, E^{\bullet} \tenp \Omega [1])$. In case $E^{\bullet}$ is
 a complex comprising a vector bundle $E$ concentrated at degree
 $0$, we check that this gives us the Atiyah class of $E$. We can
 this call this element of $\text{Hom}_{\dcat}(E^{\bullet}, E^{\bullet} \tenp \Omega
 [1])$ the Atiyah class of $E^{\bullet}$, and denote it by $\alpha(E^{\bullet})$.
 \\

 We now prove the following proposition stated in Markarian [3]
 without proof.  It has also been proven in [14].\\

 \begin{prop}

 Let $E^{\bullet}$ be a complex of vector bundles with differential $d$ on $X$ such that every
 term of $E^{\bullet}$ has a global connection $\triangledown$. Then the Atiyah
 class $\alpha(E^{\bullet})$ is given by $ \{ (-1)^n (- \triangledown d + d
 \triangledown ) \} \in \text{Hom}_{\dcat}(E^{\bullet}, E^{\bullet} \tenp \Omega
 [1]) $.

 \end{prop}

\begin{proof}
We recall that a connection $\triangledown$ on a vector bundle $E$
is a $K$-linear map \\ $\triangledown:E \rar E \otimes_{\strc}
\Omega$ so that
for a section $e$ of $E$ over an open set $U$, $\triangledown(f.e) = f(\triangledown e) + e \otimes df $. \\

Also recall that $J_1(E)= p_{2*}(p_{1}^*E \otimes {\calg O}_{X \times_K X}/{\calg I}^2) $ where $\calg I$ is the kernel of the multiplication
map $\strc \otimes_K \strc \rar \strc$. On an open subset $U= \text{Spec} R$ of $X$  where the sections of $E$ on $U$ are given by an $R$-
module $M$, $J_1(E) = M \otimes_R (R \otimes R)/I^2$ with $R$ -module structure given by multiplication with the second factor of $R \otimes_K
R$. For $a$,$b$ in $R$, let $\hat{a \otimes b}$ be the image of $a \otimes b$ in $R \otimes_K R /I^2$. Then , we have a map $p:E \rar J_1(E)$
such that $m \leadsto m \otimes \hat{1 \otimes 1}$. We observe that for a morphism of vector bundles $f: E \rar E'$ , $p
\circ f = J_1(f) \circ p$  \\

Consider a complex $E^{\bullet}$ of vector bundles on $X$ as in this
proposition.  Note that since, $E^n$ has a global connection
$\triangledown_n$ , the exact sequence $0 \rar E^n \tenp \Omega \rar
J_1(E^n) \rar E^n \rar 0$ splits. The splitting map is given by $p -
\triangledown_n $ where $p: E^n \rar J_1(E^n) $ is as in the
previous paragraph. This splitting gives rise to a map $\varphi_n:
J_1(E^n) \rar E^n \tenp \Omega$ of $\strc$-modules  We observe that
$ \{ (-1)^n \varphi_n \} $ gives us a homotopy between the second
projection \\ $J_1(E^{\bullet}) \oplus E^{\bullet} \tenp \Omega [1]
\rar E^{\bullet} \tenp \Omega [1] $ and the map  \\ $\psi :
J_1(E^{\bullet}) \oplus E^{\bullet} \tenp \Omega [1] \rar
E^{\bullet} \tenp \Omega [1]$ given by $\psi(x,y) = (-1)^{n+1}
\varphi_{n+1} dx + (-1)^n d \varphi_n x $. We next note that the map
$p: E^{\bullet} \rar J_1(E^{\bullet})$ is a map of complexes of
sheaves of $K$ vector spaces inverting the quasi-isomorphism $
\tot(0 \rar E^{\bullet} \tenp \Omega \rar J_1(E^{\bullet})) \rar
E^{\bullet}$. Therefore, in the category of complexes of sheaves of
$K$ vector spaces, $\alpha(E^{\bullet})$ is homotopy equivalent to
$\psi \circ p =  \{ (-1)^n (- \triangledown d + d \triangledown) \}
$. We note that the latter is a morphism of complexes of $\strc$
-modules. Thus $\alpha(E^{\bullet}) = \{ (-1)^n ( -\triangledown d +
d \triangledown) \} $ as a morphism in $\dcat$.

\end{proof}

We recall Kapranov [6] (Proposition 1.2.2 in 6) which amount to saying that the Atiyah class of $T_X$ which is a morphism in $\dcat$ from $T_X
[-1] \otimes T_X [-1] \rar T_X [-1] $ equips $T_X [-1]$ with the structure of a Lie algebra in $\dcat$. \\

Let $\beta$ be as in Corollary 3, Section 4. By Corollary 3 of
Section 4, \\  $\beta: T_X [-1] \rar L(\dd^1(X)) $ is a
quasiisomorphism provided that the right hand side is equipped with
the Hochschild coboundary as differential. Thus, in $\dcat$,
$T_X[-1]$ is identified with $L(\dd^1(X))$ equipped with Hochschild
co-boundary. We use
this to realize $\alpha(T_X)$ explicitly as a map of complexes from $L(\dd^1(X)) \otimes L(\dd^1(X))$ to $L(\dd^1(X))$.  \\

\begin{thm2}

The Atiyah class $\alpha(T_X)$ corresponds to the natural Lie bracket in $L(\dd^1(X))$ under the quasi-isomorphism $\beta$. In other words, the
following diagram commutes in $\dcat$ \\

$\begin{CD} T_X[-1] \tenp T_X[-1] @> \beta \otimes \beta >> L(\dd^1(X)) \tenp
  L(\dd^1(X)) \\
  @V \alpha(T_X) VV                    @VV [,] V \\
  T_X[-1] @> \beta >> L(\dd^1(X)) \\
\end{CD}$
\\

\end{thm2}

Note that all maps in the diagram in Theorem 2 except for $\alpha(TX)$ arise out of maps in $\bbcc$. The rows in the diagram of Theorem 2 are
quasi-isomorphisms. \\

  \begin{proof}
  First of all we note that $\dd^n(X)$ has a natural connection for any $n$. To see this,
 let $U = \text{Spec A}$ be an affine open subscheme of $X$. If $f \in \dd^n(U)$, then we can define $\triangledown f$ by
  $\triangledown_Y f(a_1,....,a_n) = \partial_Y f(a_1,...,a_n) \text{ } \forall \text{ } a_1,...,a_n \in A $
  It suffices to check this theorem locally. Further,
  one notes that $\alpha(T_X) = \alpha(T_X [-1]) =
  \alpha(L(\dd^1(X)))$. The second equality is by Corollary 3 Section 4. Note that
if $f \in \dd^m(U)$ and $g \in \dd^n(U)$, then $\triangledown_Y [f,g] = [\triangledown_Y f, g] +[ f, \triangledown_Y g]$.
 This implies that the connection $\triangledown$ on $\dd^n(X)$ restricts to a connection on $L(\dd^1(X)) \cap \dd^n(X)$.
  We now make the following calculation \\

  $$\tdd_Y (d \partial_I (a_1,a_2)) = \tdd_Y (a_1 \partial_I a_2 -
  \partial_I( a_1a_2) + (\partial_I a_1) a_2) $$
 $$ d (\tdd_Y (\partial_I (a_1,a_2))) = a_1 \partial_Y \partial_I a_2 -
  \partial_Y \partial_I(a_1a_2) + (\partial_Y \partial_I a_1) a_2 $$
 $$ \partial_Y (a_1 \partial_I a_2) = (\partial_Y a_1) (\partial_I
  a_2 ) +  a_1 \partial_Y \partial_I a_2 $$
 $$ \partial_Y ( (\partial_I a_1)  a_2) = (\partial_Y \partial_I a_1)
  a_2 +   (\partial_I a_1) (\partial_Y a_2) $$ \\

  This tells us that $(d \triangledown - \triangledown d) (\partial_I \otimes \partial_Y) = [\partial_I, \partial_Y] $ . \\

  Now let $D_1 \in \dd^m(U)$ and $D_2 \in \dd^n(U)$. Suppose that $(d \tdd_Y -\tdd_Y d) D_1 = {(-1)}^m [D_1, \partial_Y]$
   and $(d \tdd_Y -\tdd_Y d) D_2 = {(-1)}^n [D_2,\partial_Y]$. Then, by Proposition 2 and the
   fact that \\ $[D_1,D_2] = D_1 \otimes D_2 - {(-1)}^{mn} D_2 \otimes
D_1$,  \\

  $$ d[D_1, D_2] = [dD_1, D_2] +{(-1)}^m [D_1,dD_2] $$  Further,
  $$ \tdd_Y [D_1,D_2] = [ \tdd_Y D_1, D_2] + [D_1, \tdd_Y D_2] $$

  Thus,

  $$ (d \tdd_Y -\tdd_Y d) [D_1,D_2] = [ (d \tdd_Y -\tdd_Y d) D_1,D_2] +{(-1)}^m [D_1, (d \tdd_Y -\tdd_Y d) D_2] $$
  $$ = {(-1)}^m [[D_1, \partial_Y],D_2] + {(-1)}^{m+n} [D_1,[D_2,\partial_Y]] = {(-1)}^{m+n}[ [D_1,D_2], \partial_Y]$$

   Thus, $ {(-1)}^{m+n} (d \tdd_Y -\tdd_Y d) [D_1,D_2] = [ [D_1,D_2], \partial_Y]$. Using induction on the degree of $D \in L(\dd^1(U))$, together
   with the fact that $\dd^1(U)$ generates $L(\dd^1(U))$ as a Lie algebra over $A$, we see that ${(-1)}^{|D|} (d \tdd_Y -\tdd_Y d)D = [D,
   \partial_Y] $. This is exactly the desired theorem. \\

 \end{proof}

Corollary 1 now follows immediately from Theorem 1 and Theorem 2. \\

\section{The universal enveloping algebra of $T_X[-1]$ in $\dcat$ }

We would like to understand how Corollary 1 helps relate the big Chern classes to the Chern character. This requires Theorem 3, which is proven
in this section. A less tangible consequence of Theorem 3 is the ability to give Corollary 1 an interpretation along the lines of that given in
Section 5.2.1 to Theorem 1. \\

An associative algebra $A$ in $\dcat$ is an object of $\dcat$ such that there is a multiplication morphism $\mu_A: A \otimes A \rar A$ in
$\dcat$ which is associative i.e, \\ $\mu_A \circ (\mu_A \otimes id)= \mu_A \circ (id \otimes \mu_A)$ as morphisms in $\dcat$ from $A \otimes A \otimes A$ to $A$.\\

A Lie algebra $L$ in $\dcat$ is an object in $\dcat$ equipped with a morphism \\ $[,]_L: L \otimes L \rar L$ in $\dcat$ such that \\
(i) $[,]_L = - [,]_L \circ \tau$ where $\tau: L \otimes L \rar L \otimes L$ is the swap map. \\
(ii) $[,]_L \circ (id \otimes [,]_L) =[,]_L \circ ([,]_L \otimes id) + [,]_L \circ (id \otimes [,]_L) \circ (\tau \otimes id) $ \\

Note that any algebra $A$ in $\dcat$ has a Lie algebra structure with $[,]_A = \mu_A \circ (id -\tau)$. \\

By Proposition 2 , $\dd^{\bullet}(X)$ is an algebra in $\bbcc$ and $L(\dd^1(X))$ is a Lie algebra in $\bbcc$. Since $\dd^{\bullet}(X)$ and
$L(\dd^1(X))$ are complexes with bounded below coherent co-homology by Proposition 5 and Corollary 3 , they represent objects in $\dcat$ denoted
again by $\dd^{\bullet}(X)$ and $L(\dd^1(X))$ respectively. The algebra structure of $\dd^{\bullet}(X)$ in $\bbcc$ induces an algebra structure
in $\dcat$ . The Lie algebra structure of $L(\dd^1(X))$ in $\bbcc$ induces a Lie algebra structure in $\dcat$. Corollary 3 of Section 4 says
that $L(\dd^1(X))$ is isomorphic to $T_X[-1]$ in $\dcat$. Theorem 2 says that the Lie algebra structure on $L(\dd^1(X))$ described in this
paragraph coincides with the Lie algebra structure on $T_X[-1]$ induced by $\alpha(TX)$ after identifying $L(\dd^1(X))$ with $T_X[-1]$ via the
quasi-isomorphism $\beta$ of Corollary 3 of Section 4. \\

Given a Lie algebra $L$ in $\dcat$, its universal enveloping algebra (if it exists) is an algebra $U(L)$ in $\dcat$ together with a morphism $i:
L \rar U(L)$ of Lie algebras such that given any morphism $f: L \rar A$ of Lie algebras from $L$ to an algebra $A$ in $\dcat$, there exists a
unique morphism $\bar{f}:U(L) \rar A$ of algebras in $\dcat$ such that $f = \hat{f} \circ i$.  We now state and prove Theorem 3 \\

\begin{thm3}
$\dd^{\bullet}(X)$ is the universal enveloping algebra of $T_X[-1]$ in $\dcat$. In other words, let  $A$ be an associative algebra in $\dcat$.
If $f: T_X[-1] \rar A$ is a morphism in $\dcat$ making the following diagram commute in $\dcat$ \\

 $$\begin{CD}
 T_X[-1] \tenp T_X[-1] @> \alpha(T_X) >> T_X[-1] \\
 @V f \otimes f VV                     @VVfV \\
 A \tenp A @> [, ]_A >> A \\
 \end{CD}$$ \\

 then $\exists$ a unique morphism  $\bar{f}: \dd^{\bullet}(X) \rar A$ of algebras in
 $\dcat$ so that the composite $\bar{f} \circ \hkr = f$.
 \end{thm3}

{\it Remark :} This theorem should be compared to Lemma 2. The reason why proving this is harder is that an algebra in $\dcat$ may be realized
as a  complex of $\strc$-modules , but may not be realizable as an algebra in $\bbcc$. For instance given an algebra $A$ in $\dcat$, and a
complex $F^{\bullet}$  quasi-isomorphic to $A$, any map $F^{\bullet} \otimes F^{\bullet} \rar
F^{\bullet}$ representing $\mu_A$ may be associative only up to homotopy. \\

\subsection{Proof of Theorem 3}

Let $\beta: T_X[-1] \rar L(\dd^1(X))$ be as defined in Section 4.2. Note that $\beta$ is a quasi-isomorphism by Proposition 5, and therefore
induces an isomorphism in $\dcat$. Moreover, by Theorem 2, $\beta$ is a morphism of Lie algebras in $\dcat$. It follows that we can replace
$T_X[-1]$ by $L(\dd^1(X))$ . Suppose that $A$ is an associative algebra in $\dcat$ and that $f: L(\dd^1(X)) \rar A$ is a morphism of
Lie algebras in $\dcat$.\\

{\it Step 0:}(Conventions, Notations and some observations) \\

1. An inclusion in $\bbcc$ will mean a map of complexes that is injective term by term. \\

2. By convention, $L(\dd^1(X))$ is equipped with the Hochschild co-boundary.\\

3. If $M$ is any object in $\bbcc$, $\frel{M}$ denotes the free Lie algebra generated in $\bbcc$ over $\strc$ by $M$. The only differential on
$\frel{M}$ arises out of the differential on $M$. $\ufrel{M}$ will denote the tensor algebra generated by $M$ over $\strc$ . Again,
 the only differential on $\ufrel{M}$ arises out of the differential on $M$.   \\

4. $\sss^{\bullet}(M)$ will denote the symmetric algebra generated
by $M$ over $\strc$ in \\ $\bbcc$. We have a map of complexes
$\smap{M}: \sss^{\bullet}(M) \rar \ufrel{M}$. Given an open
subscheme $U = \text{Spec R}$ of $X$, if $m_1,...,m_k$ are sections
of $M |_U$ of degrees $d_1,...,d_k$ respectively, then
$$\smap{M}(m_1...m_k) = \frac{1}{k!} \sum_{\sigma \in S_k} s(\sigma)
m_{\sigma(1)} \otimes ... \otimes m_{\sigma(k)}$$
where $s(\sigma)$ is the sign (depending on $d_1,...,d_k$ and $\sigma$ ) such that \\ $m_1....m_k  = s(\sigma)m_{\sigma(1)}...m_{\sigma(k)} \in \sss^k(M)$. \\

Also, we have a map of complexes $\bsmap{M}:\sss^{\bullet}(\frel{M}) \rar \ufrel{M}$. If $U = \text{Spec }R$ is an open affine subscheme of $X$,
and if $z_1,...,z_k$ are sections of $\frel{M} |_U$ of degrees $d_1,...,d_k$ respectively, then $$ \bsmap{M}(z_1...z_k) = \frac{1}{k!}
\sum_{\sigma \in S_k} s(\sigma) z_{\sigma(1)} \otimes .... \otimes z_{\sigma(k)}$$ Here, the $z_i$'s on the right hand side are thought of as
sections of $\ufrel{M} |_U$ . \\ $s(\sigma)$ is the sign (depending on $d_1,...,d_k$ and $\sigma$ ) such that
\\ $z_1....z_k  = s(\sigma)z_{\sigma(1)}...z_{\sigma(k)} \in \sss^k(\frel{M})$. \\

5. Recall from Lemma 1 that the symmetrization map \\ $I: \oplus_k
\sss^k(L(\dd^1(X))) \rar \dd^{\bullet}(X)$ is an isomorphism in
$\bbcc$. Denote
its inverse by $G$. \\

6. Denote the restriction of $I$ to $L(\dd^1(X))$ by $I_1$. This is an inclusion in \\ $\bbcc$. \\

7. Denote that natural inclusion from $L(\dd^1(X))$ to $\ufrel{L(\dd^1(X))}$ by $I_2$. \\

8. The object in $\dcatb$ represented by an object $M$ in $\bbcc$ will be denoted by $M$ itself provided it exists. \\

9. The Hochschild-Kostant-Rosenberg theorem implies that all objects in \\ $\bbcc$ mentioned in this proof represent objects in $\dcatb$. \\

{\it Step 1: }(Construction of $\hat{f} : \ufrel{L(\dd^1(X))} \rar A$ in $\dcat$ ) \\

We have $f^{\otimes k}: L(\dd^1(X))^{\otimes k} \rar A^{\otimes k}$.
Since $A$ is an algebra in $\dcat$ we have the $k$ fold
multiplication $\mu_{k,A} : A^{\otimes k} \rar A$ in $\dcat$. We
thus get a morphism \\ $\mu_{k,A}  \circ f^{\otimes k}:
L(\dd^1(X))^{\otimes k} \rar A$ in $\dcat$. Putting these together
we get a morphism  $\hat{f}: \ufrel{L(\dd^1(X))} \rar A$ in $\dcat$.
By construction, $\hat{f}$ is an algebra
homomorphism in  $\dcat$ and the composite $\hat{f} \circ I_2 $ equals $f$ . \\

{\it Step 2}(Construction of $\calg I:\dd^{\bullet}(X) \rar \ufrel{L(\dd^1(X))}$

We have another kind of symmetrization map \\ $\smap{(L(\dd^1(X)))}:
\sss^{\bullet}(L(\dd^1(X))) \rar \ufrel{L(\dd^1(X))}$.  Consider the
composite \\
$\smap{L(\dd^1(X))} \circ G: \dd^{\bullet}(X) \rar \ufrel{L(\dd^1(X))}$ . Denote it by $\calg I$.  Note that $\calg I \circ I_1 = I_2$.  \\

 {\it Step 3:} (Construction of $\bar{f}: \dd^{\bullet}(X) \rar A$ in $\dcat$)\\

 We define $\bar{f} := \hat{f} \circ \calg I $.  By construction the composite
 $\bar{f} \circ I_1 = \hat{f} \circ \calg I \circ I_1 = \hat{f} \circ I_2 = f$. \\

 To complete the proof of theorem 1, we only need to check that $\bar{f}$ is an algebra homomorphism in
$\dcat$. \\

{\it Step 4 :} (Construction of an algebra homomorphism $\lambda: \ufrel{L(\dd^1(X))} \rar \dd^{\bullet}(X)$ in $\bbcc$) \\

Let $m_k:\dd^{\bullet}(X)^{\otimes k} \rar \dd^{\bullet}(X)$ denote the $k$-fold product of $\dd^{\bullet}(X)$. \\

Consider the composite map

$$\lambda_k: \begin{CD} L(\dd^1(X))^{\otimes k} @> {I_1}^{\otimes k}>> \dd^{\bullet}(X)^{\otimes k} @> m_k >> \dd^{\bullet}(X) \end{CD} $$  This map is a map
in $\bbcc$ as $I_1$ and $m_k$ are maps in $\bbcc$. \\
 Putting these together, we obtain a map $\lambda: \ufrel{L(\dd^1(X))} \rar \dd^{\bullet}(X)$. By construction, $\lambda$ is a morphism of algebras in $\bbcc$.
  Also, the composite of $\lambda_k$ with the restriction of $\smap{L(\dd^1(X))}$ to $\sss^k(L(\dd^1(X)))$ is just the restriction of the
  symmetrization map $I$ to $\sss^k(L(\dd^1(X)))$ Thus, $\lambda \circ \smap{L(\dd^1(X))} = I$. It follows that
  \\ $\lambda \circ \calg I = \lambda \circ \smap{L(\dd^1(X))} \circ G = I \circ G =id$. \\

{\it Step 5:} (Uniqueness of $\bar{f}$ ) \\

Suppose that $\bar{f}_1:\dd^{\bullet}(X) \rar A$ and
$\bar{f}_2:\dd^{\bullet}(X) \rar A$ are two morphisms of algebras in
$\dcat$ such that $\bar{f}_1 \circ I_1 = \bar{f}_2 \circ I_1 = f$.
Then, by the construction of the map $\lambda:\ufrel{L(\dd^1(X))}
\rar \dd^{\bullet}(X)$ in Step 4 ,  $$\bar{f}_1 \circ \lambda =
\bar{f}_2 \circ \lambda = \hat{f} :\ufrel{L(\dd^1(X))} \rar A \text{
. }$$ It follows from this and the fact (demonstrated in Step 4)
that $\lambda \circ \calg I = id$ that
$$\hat{f} \circ \calg I = \bar{f}_1 = \bar{f}_1 \circ \lambda \circ \calg I =
\bar{f}_2 \circ \lambda \circ \calg I = \bar{f}_2 \text{ . }$$ This
proves that $\bar{f}_1$ and $\bar{f}_2$ are identical morphisms in
$\dcat$.\\

{\it Step 6: } We now return to proving that $\bar{f}$ is a morphism
of algebras in $\dcat$. Proving this will completely prove Theorem
3. We claim that the following diagram commutes in $\bbcc$ \\

$$ \begin{CD} \dd^{\bullet}(X) \otimes \dd^{\bullet}(X) @> \calg I \otimes \calg I >>
\ufrel{L(\dd^1(X))} \otimes \ufrel{L(\dd^1(X))} \\
@V m VV                        @VV \mu_{\ufrel{L(\dd^1(X))}} V \\
\dd^{\bullet}(X) @<< \lambda <    \ufrel{L(\dd^1(X))}\\
\end{CD} $$ \\

This is immediate from the fact that $\lambda \circ \calg I = id$ and the fact that $\lambda$ is an algebra homomorphism in $\bbcc$.
 Both these facts were demonstrated in Step 4. \\

It follows that the following diagram commutes in $\dcat$ \\

$$ \begin{CD} \dd^{\bullet}(X) \otimes \dd^{\bullet}(X) @> \calg I \otimes \calg I >> \ufrel{L(\dd^1(X))} \otimes \ufrel{L(\dd^1(X))}  @>
\hat{f} \otimes \hat{f}>> A
\otimes A \\
@V m VV       @V \mu_{\ufrel{L(\dd^1(X))}} VV    @VV \mu_A V \\
\dd^{\bullet}(X) @<< \lambda <  \ufrel{L(\dd^1(X))}  @> \hat{f} >> A \\
\end{CD} $$ \\

The arrows in the square on the left arise out of morphisms in $\bbcc$ itself.\\

Now observe that in $\dcat$ , $$ \mu_A \circ (\bar{f} \otimes \bar{f}) = \mu_A \circ [(\hat{f} \circ \calg I) \otimes (\hat{f} \circ \calg I)]
$$
$$ = \hat{f} \circ \mu_{\ufrel{L(\dd^1(X))}} \circ (\calg I \otimes \calg I) = \hat{f} \circ (id - \calg I \circ \lambda) \circ \mu_{\ufrel{L(\dd^1(X))}}
 \circ (\calg I \otimes \calg I) $$
$$ + \hat{f} \circ \calg I \circ \lambda \circ \mu_{\ufrel{L(\dd^1(X))}} \circ (\calg I \otimes \calg I) $$
But $\lambda \circ \mu_{\ufrel{L(\dd^1(X))}} \circ (\calg I \otimes \calg I) = m $ . Thus, \\
$$\mu \circ \bar{f} \otimes \bar{f} = \hat{f} \circ (id - \calg I \circ \lambda) \circ \mu_{\ufrel{L(\dd^1(X))}} \circ (\calg I \otimes \calg I)+ \hat{f} \circ \calg I \circ
m $$ $$ =  \hat{f} \circ (id - \calg I \circ \lambda) \circ \mu_{\ufrel{L(\dd^1(X))}} \circ (\calg I \otimes \calg I)+ \bar{f} \circ m $$

 Therefore, to show that $\bar{f}$ is a homomorphism of algebras in $\dcat$, it suffices to show the following proposition. \\

\begin{prop}
The composite  $$\begin{CD} \ufrel{L(\dd^1(X))} @> id - \calg I \circ \lambda
>> \ufrel{L(\dd^1(X))} @> \hat{f} >> A \end{CD} $$  is $0$ in $\dcat$. \\
\end{prop}

We note that the first arrow $id - \calg I \circ \lambda$ is an arrow arising out of a morphism in $\bbcc$. We break the proof
of Proposition 7 down into the following easier propositions \\

Note that any morphism in $\bbcc$ can be thought of as morphism of graded $\strc$ modules.
Denote the map $\bsmap{L(\dd^1(X))} :\sss^{\bullet}(\frel{L(\dd^1(X))}) \rar \ufrel{L(\dd^1(X))}$ by $B$. Note
that $B$ is a morphism in $\bbcc$. Denote the restriction of $B$ to $\sss^{k}(\frel{L(\dd^1(X))})$ by $B^k$. In particular,
$B^1$ is a map in $\bbcc$ from $\frel{L(\dd^1(X))}$ to $\ufrel{L(\dd^1(X))}$. \\

\begin{prop} $B^1$ is a morphism of Lie algebras in $\bbcc$ \end{prop}

\begin{proof} Let $\bsmap{M}^k$ denote the restriction of $\bsmap{M}$ to $\sss^k(\frel{M})$. Then , we claim that $\bsmap{M}^1$ is a morphism
of Lie algebras in $\bbcc$. Let $U= \text{Spec R}$ be an open affine subscheme of $X$ and let $z_1,z_2$ be sections of $\frel{M} |_U$ of degrees
$d_1,d_2$ respectively. Then $\bsmap{M}^1[z_1,z_2] = z_1 \otimes z_2 - {(-1)}^{d_1d_2} z_2 \otimes z_1 $ and $\bsmap{M}^1(z_i) = z_i$. This
proves the desired proposition.
\end{proof}

Recall that we have a natural inclusion of complexes $L(\dd^1(X)) \rar \frel{L(\dd^1(X))}$. This just treats a section of $L(\dd^1(X))$ as a
section of $\frel{L(\dd^1(X))}$. \\

\begin{prop} (i) As a morphism of graded $\strc$ modules, $\lambda \circ B^1 $ maps $\frel{L(\dd^1(X))}$ to $I_1(L(\dd^1(X)))$.
Let $\pi$ denote $G \circ \lambda \circ B^1$.\\
(ii) $\pi: \frel{L(\dd^1(X))} \rar L(\dd^1(X))$ is a map of Lie algebras in $\bbcc$. \\
(iii) The composite $ \begin{CD} L(\dd^1(X)) @>>> \frel{L(\dd^1(X))} @>{\pi}>> L(\dd^1(X)) \end{CD} $ is the identity. \\
(iv) $\lambda \circ B = I \circ  \sss^{\bullet}(\pi)$ \\
\end{prop}

\begin{proof} Note that $\ufrel{L(\dd^1(X))}$ and $\dd^{\bullet}(X)$ are  Lie algebras in $\bbcc$. Their Lie algebra structures
are induced by their algebra structures in $\bbcc$ . Since $\lambda: \ufrel{L(\dd^1(X))} \rar \dd^{\bullet}(X)$ is a morphism of algebras in
$\bbcc$, $\lambda$ is also a morphism of Lie algebras in $\bbcc$. Also, by its very construction, $\lambda$ maps $I_2(L(\dd^1(X)))$ identically
to $I_1(L(\dd^1(X)))$. It follows that $\lambda$ maps $B^1(\frel{L(\dd^1(X))})$ to the Lie subalgebra of
$\dd^{\bullet}(X)$ generated over $\strc$ by $I_1(L(\dd^1(X)))$ which is  $I_1(L(\dd^1(X)))$ itself. This proves (i)\\

Observe that $\lambda$ is a Lie algebra homomorphism in $\bbcc$.  Further, since $I_1:L(\dd^1(X)) \rar I(L(\dd^1(X)))$ is a morphism of Lie
algebras in $\bbcc$ ,  $G:I(L(\dd^1(X))) \rar L(\dd^1(X))$ is also a morphism of Lie algebras in $\bbcc$. Since $B^1$ is  a morphism of Lie
algebras in $\bbcc$ by Proposition 8, (ii)
follows. \\

(iii) is immediate from the construction of $\lambda$ . \\

Let $U= \text{Spec } R$ be an affine open subscheme of $X$. Let $x_1,....,x_k$ be homogenous sections of $\frel{L(\dd^1(U))}$ of degrees
$d_1,....,d_k$ respectively . For a permutation $\sigma \in S_k$ let $s(\sigma)$ be the sign such that $x_1...x_k = s(\sigma)
x_{\sigma(1)}....x_{\sigma(k)}$ in $\sss^k(\frel{L(\dd^1(U))})$. Then,
$$B(x_1...x_k ) = \frac{1}{k!} \sum_{\sigma \in S_k} s(\sigma) x_{\sigma(1)} \otimes ... \otimes x_{\sigma(k)}$$ $$\lambda \circ B(x_1....x_k) = \frac{1}{k!}
 \sum_{\sigma \in
S_k} s(\sigma) \lambda(x_{\sigma(1)}) \otimes .... \otimes \lambda(x_{\sigma(k)})$$. But by the definition of $\pi$, if $x \in
\frel{L(\dd^1(U))}$, then $\lambda(x) = I_1(\pi(x))$. It follows that $\lambda \circ B (x_1...x_k)$ is precisely  $I(\pi(x_1)....\pi(x_k))$ in
$\dd^{\bullet}(X)$. This proves (iv).

\end{proof}

Recall that we have the map \\ $\smap{\frel{L(\dd^1(X))}}:
\sss^{\bullet}(\frel{L(\dd^1(X))}) \rar \ufrel{\frel{L(\dd^1(X))}}$.
Let $\smap{\frel{L(\dd^1(X))}}^k$ denote its restriction to
$\sss^k(\frel{L(\dd^1(X))})$. Let $\smap{L(\dd^1(X))}^k$ denote the
restriction of $\smap{L(\dd^1(X))}$ to $\sss^k(L(\dd^1(X)))$.  \\

Let $\pi^{\bullet}$ denote the morphism $\oplus_k \pi^{\otimes k}: \ufrel{\frel{L(\dd^1(X))}} \rar \ufrel{L(\dd^1(X))}$ \\

\begin{prop}
The following diagram commutes in $\bbcc$  \\
$$\begin{CD}
{\text{Sym}}^k(\frel{L(\dd^1(X))}) @> \smap{\frel{L(\dd^1(X))}}^k >> \ufrel{\frel{L(\dd^1(X))}} \\
@VV{{\text{Sym}}^k(\pi)}V @VV{\pi^{\bullet}}V \\
\sss^k(L(\dd^1(X))) @> \smap{L(\dd^1(X))}^k  >> \ufrel{L(\dd^1(X))} \\
\end{CD} $$
\end{prop}

\begin{proof}
This is immediate from the definitions of $\pi$, $\smap{\frel{L(\dd^1(X))}}^k$ and $\smap{L(\dd^1(X))}^k$. \\
\end{proof}

Note that the symmetrization map \\ $\bsmap{L(\dd^1(X))}:
\sss^{\bullet}(\frel{L(\dd^1(X))}) \rar \ufrel{L(\dd^1(X))}$ is an
isomorphism in \\ $\bbcc$. This is proven in a fashion similar to
Lemma 1. Denote the inverse of $\bsmap{L(\dd^1(X))}$ by $W$. Let
$\calg Z$
denote $\smap{\frel{L(\dd^1(X))}} \circ W$. We have the following corollary of Propositions 9 and 10. \\

\begin{cor}

The following diagram commutes in the category of differential graded $\strc$ modules \\
$$ \begin{CD}
\ufrel{L(\dd^1(X))} @> \calg Z >> \ufrel{\frel{L(\dd^1(X))}} \\
@VV{\lambda}V               @VV{\pi^{\bullet}}V \\
\dd^{\bullet}(X)    @> \calg I >> \ufrel{L(\dd^1(X))} \\
\end{CD} $$
\end{cor}

\begin{proof}
Observe that part (iv) of  Proposition 9 tells us that $\lambda
\circ B = I \circ \sss^{\bullet}(\pi)$. Thus, $G \circ \lambda \circ
B = \sss^{\bullet}(\pi)$. Thus, $$ \smap{L(\dd^1(X))} \circ G \circ
\lambda \circ B \circ W = \smap{L(\dd^1(X))} \circ
\sss^{\bullet}(\pi) \circ W$$ But $\calg I = \smap{L(\dd^1(X))}
\circ G$ and $B \circ W = id$. Thus, $\calg I \circ \lambda =
\smap{L(\dd^1(X))} \circ \sss^{\bullet}(\pi) \circ W$. By
Proposition 10, $\smap{L(\dd^1(X))} \circ \sss^{\bullet}(\pi) =
\pi^{\bullet} \circ \smap{\frel{L(\dd^1(X))}}$. Thus,\\ $\calg I
\circ \lambda = \pi^{\bullet} \circ \smap{\frel{L(\dd^1(X))}} \circ
W = \pi^{\bullet} \circ \calg Z$.
\end{proof}

We have the multiplication map $\mu_{\calg T}$ from
$\ufrel{\ufrel{L(\dd^1(X))}}$ to $\ufrel{L(\dd^1(X))}$. This arises
out of the tensor product in $\ufrel{L(\dd^1(X))}$. Let \\ $B:
\sss^{\bullet}(\frel{L(\dd^1(X))}) \rar \ufrel{L(\dd^1(X))}$ be as
in Proposition 9.  The map $B^1$ yields us a map $\ufrel{B^1}:
\ufrel{\frel{L(\dd^1(X))}} \rar \ufrel{\ufrel{L(\dd^1(X))}}$ in
$\bbcc$. This is a map of algebras in $\bbcc$ by construction.
Denote the composite $\mu_{\calg T} \circ \ufrel{B^1} :
\ufrel{\frel{L(\dd^1(X))}}
\rar \ufrel{L(\dd^1(X))}$ by $\calg G$. Since $\mu_{\calg T}$ and $\ufrel{B^1}$ are morphisms of algebras in $\bbcc$ so is $\calg G$. \\

\begin{prop} With $\calg G$ as defined above,  $\calg G \circ \calg Z = id$. \\
\end{prop}

\begin{proof}
 First note that $\calg G \circ \smap{\frel{L(\dd^1(X))}} = \bsmap{\frel{L(\dd^1(X))}}$. To see this, let \\ $U = \text{Spec R}$ be an
 open affine subscheme of $X$, and let $z_1,...,z_k$ be sections of \\ $\frel{L(\dd^1(X))} |_U$ of degrees $d_1,...,d_k$ respectively. Then, by the
 construction of $G$,\\ $G(z_1 \otimes ... \otimes z_k) = z_1 \otimes ... \otimes z_k$ where the tensor product on the right is that in
 \\ $\ufrel{L(\dd^1(X))}$ and where the $z_i$'s on the right are treated as sections of $\ufrel{L(\dd^1(X))} |_U$.

Now, $\calg G \circ \calg Z = \calg G \circ \smap{\frel{L(\dd^1(X))}} \circ W = \bsmap{\frel{L(\dd^1(X))}} \circ W = id $ .\\

\end{proof}

The following corollary is obtained from the above proposition and Corollary 4.

\begin{cor} $id - \calg I \circ \lambda = (\calg G - \pi^{\bullet}) \circ \calg Z $
\end{cor}

Note that all the commutative diagrams in Propositions 8-11 and their corollaries are diagrams in $\bbcc$. They induce corresponding commutative
diagrams in $\dcat$. \\

From corollary 5, it is clear that to prove Proposition 7, it suffices to prove the following proposition. \\

\begin{prop} $\hat{f} \circ \calg G =  \hat{f} \circ \pi^{\bullet} $ in $\dcat$.
\end{prop}

\begin{proof} We have a natural inclusion of complexes $\frel{L(\dd^1(X))}^{\otimes k} \rar \ufrel{\frel{L(\dd^1(X))}}$. Denote the composite of
$\calg G$ with this inclusion by ${\calg G}_k$. This is a morphism in \\ $\bbcc$. \\

To prove this proposition ,it suffices to show that $\hat{f} \circ {\calg G}_k = \hat{f} \circ \pi^{\otimes k}$ in $\dcat$. Since $\calg G$ is a
homomorphism of algebras in $\bbcc$, ${\calg G}_k = {\calg G}_1^{\otimes k}$. Moreover, $\hat{f}$ is a morphism of algebras in $\dcat$. It is
therefore sufficient to show that $\hat{f} \circ {\calg G}_1 = \hat{f} \circ \pi$. \\

This is done in Propositions 13 and 14 that follow. \\

\end{proof}

We recall that if $V$ is a vector space over a field of characteristic $0$, $S_n$ acts on the right on $V^{\otimes n}$. If $\sigma$ is a
permutation on $S_n$, $\sigma(v_1 \otimes ... \otimes v_n) = v_{\sigma(1)} \otimes... \otimes v_{\sigma(n)}$. This extends to an action of $S_n$
on ${\calg A}^{\otimes n}$ for any complex $\calg A$ of $\strc$ modules. The action descends to an action on ${\calg B}^{\otimes n}$ for any
element $\calg B$ of $\dcat$. \\

Observation 1. Note that if $A$ is an associative algebra in $\dcat$, with multiplication $\mu_A$, the Lie Bracket $[,]:A \otimes A \rar A$ is
defined as $\mu_A \circ(id - ( 1 2))$ where $( 1 2)$ is the swap applied to $A^{\otimes 2}$. Let $\tau_k \in S_n$ be the $k$-cycle $(n-k+1
\text{ } n-k+2 \text{ } .... n)$. Let $e_n$ be the element $(id - \tau_n) \circ .... \circ (id -\tau_2) \in KS_n$. $e_n$ is a quasi-idempotent,
in the group ring $KS_n$  of $S_n$. If $\mu_{n,A}:A^{\otimes n} \rar A$ denotes the $n$-fold multiplication on $A$, then $\mu_{n,A} \circ e_n:
A^{\otimes n} \rar A$
is the $n$-fold Lie bracket on $A$. Note that this is only a morphism in $\dcat$.  \\

Denote by $L_n$ the morphism in $\bbcc$ from $L(\dd^1(X))^{\otimes n}$ to $L(\dd^1(X))$ such that for sections $z_1,...,z_n$ of $L(\dd^1(X))$
over an affine open subscheme $U$ of $X$, $L_n(z_1 \otimes ... \otimes z_n) = [z_1,[z_2,[...,[z_{n-1},z_n]]]]$. \\

The following proposition is a direct consequence of the fact that $f$ is a morphism of Lie algebras in $\dcat$ \\

\begin{prop}
The following diagram commutes in $\dcat$ \\

$$\begin{CD}
L(\dd^1(X))^{\otimes n} @>{f}^{\otimes n}>> A^{\otimes n} \\
@VV{L_n}V             @VV{\mu_{n,A} \circ e_n}V \\
L(\dd^1(X)) @>f>> A \\
\end{CD}$$

\end{prop}

\begin{proof} For $n=2$ this is exactly the statement that $f$ is a morphism of Lie algebras in $\dcat$. For other $n$ it is proven by induction
on $n$ using the facts that $L_2 \circ (id \otimes L_{n-1}) = L_n$ and $\mu_{n,A} \circ e_n = (\mu_A \circ e_2) \circ (id \otimes [\mu_{n-1,A} \circ e_{n-1}])$  \\
\end{proof}

\begin{prop}
The following diagram also commutes in $\bbcc$, and hence in $\dcat$ \\
$$\begin{CD}
\frel{L(\dd^1(X))} @>B^1>> \ufrel{L(\dd^1(X))} \\
@VV{\pi}V       @VV{\oplus_n \frac{1}{n} L_n}V \\
 L(\dd^1(X))  @>{id}>> L(\dd^1(X)) \\
 \end{CD}
 $$
\end{prop}

\begin{proof}

Let $L$ denote $\oplus_n \frac{1}{n} L_n$. Denote the multiplication on $\ufrel{L(\dd^1(X))}$ by $m_{\calg T}$. The $n$ -fold multiplication
will be denoted by $m_{n, \calg T}: \ufrel{L(\dd^1(X))}^{\otimes n} \rar \ufrel{L(\dd^1(X))}$. All these are maps in $\bbcc$. \\

Let $e_n$ be as in Observation 1 prior to Proposition 13. It is immediate from Theorem 8.16 of Reutenauer [11] that $\frac{1}{n} e_n$ is an
idempotent. Moreover, it is a projection from $\ufrel{W}$ to $\frel{W}$ for any $W \in \bbcc$. Let $z_1,...,z_n$ be sections of $L(\dd^1(X))$
over an open affine subscheme $U$ of $X$. Then, $B^1([z_1,[z_2[....[z_{n-1},z_n]]]]) = e_n(z_1 \otimes ... \otimes z_n)$ by the definition of
$e_n$. On the other hand, $\pi(B^1([z_1,[z_2[....[z_{n-1},z_n]]]])) = [z_1,[z_2[....[z_{n-1},z_n]]]]$ where the bracket on the right hand side
is that of $L(\dd^1(X))$. To verify this proposition, it suffices to check that $L_n \circ \frac{1}{n} e_n = L_n$. Let
$C:B^1(\frel{L(\dd^1(X))}) \rar \frel{L(\dd^1(X))}$ be the left inverse of $B^1$. Now, $L_n(z_1 \otimes ... \otimes z_n) = \pi \circ C \circ
e_n(z_1 \otimes.... \otimes z_n)$. Thus, $L_n \circ \frac{1}{n} e_n = \pi \circ C \circ e_n \circ \frac{1}{n} e_n = \pi \circ C \circ e_n =
L_n$.
This proves the desired proposition. \\

\end{proof}

\begin{proof} (Final steps to proving Proposition 7) \\

Note that $f \circ \pi = \hat{f} \circ \pi$. Now combining Propositions 13 and 14, we see that $\hat{f} \circ \pi$
equals the following composition of morphisms in $\dcat$ \\

$$\begin{CD} \frel{L(\dd^1(X))} @>B^1>> \ufrel{L(\dd^1(X))} @>{\oplus_n f^{\otimes n}}>> T(A) @>{\oplus_n (\mu_{n,A} \circ e_n) }>> A \end{CD}$$

Now we can see that if $h:\calg B \rar \calg C$ is any morphism in $\dcat$ and $\sigma \in S_n$ then $h^{\otimes n} \circ \sigma = \sigma \circ
h^{\otimes n}$. This is verified by checking the corresponding fact at the level of complexes of $\strc$ modules. \\

Let $\mu_{\bullet,A} = \oplus_n \mu_{n,A} :T(A) \rar A$. It follows that $f \circ \pi$ is given by the following composition of morphisms in $\dcat$ \\

$$ \begin{CD} \frel{L(\dd^1(X))} @>B^1>> \ufrel{L(\dd^1(X))} @>{\oplus_n \frac{1}{n} e_n}>> \ufrel{L(\dd^1(X))} @>{\oplus_n f^{\otimes n}}>> T(A)  @>>{\mu_{\bullet, A}}>  A
\end{CD} $$

But ${\oplus_n \frac{1}{n} e_n }$ is a projector from $T(W)$ to $L(W)$ for any complex of $\strc$ modules $W$. This is immediate from Theorem
8.16 of Reutenauer [11]. It follows that the composite of the first two maps in the previous composition is just ${\calg G}_1$ .  The map
$\mu_{\bullet, A} \circ \oplus_n f^{\otimes n}$ is precisely $\hat{f}$.
This proves that $\hat{f} \circ {\calg G}_1 = \hat{f} \circ \pi$ thereby finally proving Proposition 7 and therefore, Theorem 1. \\

\end{proof}

Remark: Taking $\calg I \circ \lambda$ essentially amounts to taking brackets among elements of $L(\dd^1(X))$ . The fact that $\bar{f}$ when
restricted to $L(\dd^1(X))$ is a Lie algebra homomorphism implies that $\bar{f} \circ \calg I \circ \lambda = \bar{f} $ in $\dcat$ which is
exactly what we want. This is the "hand waving" argument for Proposition 7 that is made rigorous by the  proofs of Propositions 8-14. \\

\section{The Chern character as a character of a representation}

If $E$ is a vector bundle on $X$, $\ennnn(E)$ concentrated in degree $0$ is an algebra in $\bbcc$. The algebra structure on $\ennnn(E)$ is given
by the composition map $\circ: \ennnn(E) \otimes \ennnn(E) \rar \ennnn(E)$. We recall from Kapranov [6] that the Atiyah class of $E$ endows it
with the structure of a module over the Lie algebra $T_X[-1]$ in $\dcat$. In other words, if we identify $\text{Hom}_{\dcat}(E \otimes T_X[-1],
E)$
with $\text{Hom}_{\dcat}(T_X[-1], \ennnn(E))$, then $\alpha(E): T_X[-1] \rar \ennnn(E)$ is a morphism of Lie algebras in $\dcat$. \\

The multiplication on $\ennnn(E)$ induces a $k$-fold multiplication $\circ_k  : \ennnn(E)^{\otimes k} \rar \ennnn(E)$. For an element $\alpha
\in \text{Hom}_{\dcat}(T_X[-1], \ennnn(E))$, let $\alpha^{\circ k}$ denote the composite
$$\begin{CD} T_X^{\otimes k}[-k] @>{\alpha^{\otimes k}}>> \ennnn(E)^{\otimes k} @>{\circ_k}>> \ennnn(E) \end{CD}$$ which is an element of
$\text{Hom}_{\dcat}(T_X^{\otimes k}[-k], \ennnn(E))$. Denote $\alpha(E)^{\circ k}$ by $\tilde{\ttt_k}(E)$. \\

Let $p: \texal{X}{k} \rar T_X^{\otimes k}[-k]$ be the morphism of
complexes such that if $v_1,...,v_k$ are sections of $T_X$ over an
open affine subscheme $U$ of $X$, then \\ $p(v_1 \wedge ... \wedge
v_k) = \frac{1}{k!} \sum_{\sigma \in S_k} \text{sgn}(\sigma)
v_{\sigma(1)} \otimes ... \otimes v_{\sigma(k)} $.\\ Then,
$\alpha(E)^{\circ k} \circ p \in \text{Hom}_{\dcat}(\texal{X}{k},
\ennnn(E))$. Denote $\alpha(E)^{\circ k} \circ p$ by
$\tilde{\chn_k}(E)$. \\

Note that we have a map of $\strc$-modules $tr: \ennnn(E) \rar \strc$. This is a map in $\bbcc$ if $\ennnn(E)$ and $\strc$ are thought of as
complexes concentrated in degree $0$. Then $tr \circ \tilde{\ttt_k}(E)$ is an element in $\text{Hom}_{\dcat}(T_X^{\otimes k}[-k], \strc)$.
Denote it by $\ttt_k(E)$. Similarly, $tr \circ \chn_k(E)$ is an element of $\text{Hom}_{\dcat}(\texal{X}{k}, \strc)$. Denote it by $\chn_k(E)$.
The isomorphism of $\text{Hom}_{\dcat}(\texal{X}{k}, \strc)$ with $\text{H}^k(X,\wedge^k \Omega)$ maps $\chn_k(E)$ to the degree $k$ component
of the Chern character of $E$. \\

 Let $\hkr^1:T_X[-1] \rar \dd^{\bullet}(X)$ denote the composite of $\hkr$ with the inclusion of $T_X[-1]$ in $\oplus_k \texal{X}{k}$ as a
direct summand.Since $\alpha(E): T_X[-1] \rar \ennnn(E)$ is a morphism of Lie algebras in $\dcat$, Theorem 3 implies that there exists a
morphism
 $\theta_E:\dd^{\bullet}(X) \rar \ennnn(E)$ of algebras in $\dcat$ such that $\theta_E \circ \hkr^1 = \alpha_E$. Let
 $\varphi_E := tr \circ \theta_E:\dd^{\bullet}(X) \rar \strc$. \\

 Let $\bhkr:\oplus_k T_X^{\otimes k}[-k] \rar \dd^{\bullet}(X)$ be the morphism of complexes such that if $v_1,...,v_k$ are sections of $T_X$ over an open affine
subscheme $U$ of $X$, then \\ $\bhkr(v_1 \otimes .... \otimes v_k) = v_1 \otimes .... \otimes v_k \in \dd^n(U) $. \\

We now have the following theorem. \\

\begin{thm4}

(i) $$\tilde{\ttt_k(E)} = \theta_E \circ \bhkr $$ (ii) $$ \tilde{\chn_k(E)} = \theta_E \circ \hkr  $$ (iii) $$ \oplus_n \ttt_n(E) = \varphi_E
\circ \bhkr $$ (iv) $$ \chn(E) = \varphi_E \circ \hkr $$

\end{thm4}

\begin{proof}
Let $\bhkr^n$ denote the composite of $\bhkr$ with the inclusion of $T_X^{\otimes n}[-n]$ in $\oplus_k T_X^{\otimes k}[-k]$ as a direct summand.
Then, by the definition of $\bhkr$, $\bhkr^n = {\hkr^1}^{\otimes n}$. Further, $\theta_E$ is a morphism of algebras
in $\dcat$ and $\theta_E \circ \hkr^1 = \alpha(E)$. It follows that the following diagram commutes in $\dcat$. \\

$$ \begin{CD} T_X[-1]^{\otimes n} @> {\bhkr^n}>> \dd^{\bullet}(X) \\
@V id VV      @VV \theta_E V \\
T_X[-1])^{\otimes n} @> \alpha(E)^{\circ n} >> \ennnn(E) \\
\end{CD}$$ \\

Also, the following diagram commutes in $\dcat$ by the definitions of $\hkr$ and $\bhkr$ \\

$$ \begin{CD} \texal{X}{n} @> p >> T_X[-1]^{\otimes n} \\
@V \hkr VV                    @VV {\bhkr_n}V \\
\dd^{\bullet}(X) @< {id} << \dd^{\bullet}(X))  \\
\end{CD}$$ \\

It follows that $\theta_E \circ \bhkr =  \oplus_k \alpha(E)^{\circ k} \circ id $  \\
Thus, $\varphi_E \circ \bhkr = \oplus_k tr \circ \alpha(E)^{\circ k}) = \oplus_k \ttt_k(E)$ \\

This proves (i) and (iii). For (ii) and (iv) , we use the fact that $\bhkr \circ p = \hkr $ to see that $\theta_E \circ \bhkr \circ p = \oplus_k
\alpha(E)^{\circ k} \circ id \circ p$. Composing both sides of this with the trace map from $\enn(E)$ to $\strc$, we see that
$\varphi_E \circ \hkr = tr \circ \oplus_k \alpha(E)^{\circ k} \circ p$. This proves (ii) and (iv). \\

\end{proof}

In the classical situation, if $\mathfrak{g}$ is an ordinary Lie algebra over a field of characteristic $0$, and $E$ is a finite dimensional
representation of $\mathfrak{g}$, we have a Lie algebra homomorphism $\theta_E: \mathfrak{g} \rar \enn(E)$. This induces a map $U{\mathfrak{g}}
\rar \enn(E)$ of algebras where $U{\mathfrak{g}}$ is the universal enveloping algebra of $E$. One has the trace $\enn \rar K$. One can therefore
compose these to get a map $\varphi_E: U{\mathfrak{g}} \rar K$. This is the character of the representation $E$ of $\mathfrak{g}$. The analogy
with the Chern character is now clear. By Kapranov [6]  any vector bundle $E$ is a representation of the Lie algebra $T_X[-1]$ in $\dcat$. By
part (iv) of Theorem 4, the Chern character of the vector bundle is the character in the Representation theoretic sense of the representation
$E$ of $T_X[-1]$. \\

Theorem 4 also enables us to prove some properties of the big Chern classes shown by Ramadoss [8] in a more general framework . In this section
we shall reprove the fact that the big Chern classes commute with the Adams operations. The new proof would make the parallel of this fact with
the Representation theoretic identity $\chi_{{\psi}^p V}(g) = \chi_V
(g^p)$ transparent. \\

For this we need a digression on Adams operations in commutative
Hopf algebras. Let $\mathcal{H}$ be a commutative Hopf algebra, with
multiplication $\mu$ and comultiplication $\Delta$. Let $\mu_k$ and
$\Delta^k$ denote the $k$ fold multiplication and $k$ fold
comultiplication respectively. Then, the maps $\psi^k : \mathcal{H}
\rar \mathcal{H}$ are ring homomorphisms. Moreover, $\psi^p \circ
\psi^q = \psi^{pq}$. These maps can therefore be thought of a Adams
operations. If $\mathcal{H}$ is noncommutative but is cocommutative
and primitively generated then these maps are not ring
homomorphisms though they satisfy $\psi^p \circ \psi^q = \psi^{pq}$. \\

Coming back to the classical picture recall that the character of a representation $E$ of a Lie algebra $\mathfrak{g}$ is a $K$ linear map
$\varphi_E: U({\mathfrak{g}}) \rar K$. Let us restrict $\varphi_E$ to elements of $U{\mathfrak{g}}$ that are of the form $exp(tv)$ , $t \in
\mathbb{R}$, $ v \in \mathfrak{g}$ ,$v$ fixed. This  yields a character (in the usual "character of a representation of a group" sense) of the
representation $E$ of the one parameter group $\exx(tv) \subset U{\mathfrak{g}}$. In this case we denote the character $\varphi_E$ restricted to
the one parameter group by $\chi_E$ to keep notation more standard. \\

 From the fact that  $\chi_{{\psi}^p E}(g) =
 \chi_E(g^p)$  it follows that  $\chi_{{\psi}^p E} \exx(tv) = \chi_E \exx(ptv)$ . If $*$ denotes the multiplication in $U{\mathfrak{g}}$ and
 $v^k := v*...*v \text{ k times }$ then \\ $\chi_{\psi^p E} \sum_k \frac{1}{k!} t^k v^k = \sum_k \frac{1}{k!} t^k p^k v^k$. Note that $U{\mathfrak{g}}$
 is a co-commutative Hopf algebra with \\ $\Delta^p(v) = v \otimes 1 \otimes ... \otimes 1 + .... + 1 \otimes ... \otimes 1 \otimes v \in U{\mathfrak{g}}^{\otimes p} $ for $v \in \pbw(\mathfrak{g})$.
 Since $\Delta^p$ is an algebra homomorphism, it follows that
 $$\Delta^p(v^k) = \sum_{\{ (k_1,...,k_p) | k_i \geq 0 \text{ } \forall \text{ } i \text{and} \sum k_i=k \}} v^{k_1} \otimes ... \otimes v^{k_p}
  \in   U{\mathfrak{g}}^{\otimes p} $$ Thus, $\mu_p \circ \Delta^p (v^k) = p^k v^k$. Since $\{ \pbw(v^k) | v \in \mathfrak{g}, k \geq 0\}$ spans
  $U{\mathfrak{g}}$, the following diagram commutes   \\

 $$ \begin{CD}
  U{\mathfrak{g}} @> \psi^p >> U{\mathfrak{g}} \\
  @VidVV                     @VV\varphi_EV \\
  U{\mathfrak{g}} @> \varphi_{\psi^p E} >> K \\
  \end{CD} $$ \\

  The commuting of this diagram is equivalent to the fact that $\chi_{{\psi}^p E}(g) =
 \chi_E(g^p)$. We will now prove the same in $\dcat$ for the Lie
 algebra $T_X[-1]$. Recall from Proposition 2 that $\dd^{\bullet}(X)$ is a Hopf-algebra in $\bbcc$ . It
 is therefore, a Hopf algebra in $\dcat$ as well, with all operations induced by the corresponding operations in $\bbcc$. It is co-commutative but
 non-commutative. \\

 \begin{prop}
 The following diagram commutes in $\dcat$ \\
 $$ \begin{CD}
   \dd^{\bullet}(X) @> \psi^p>>  \dd^{\bullet}(X) \\
   @V\varphi_{\psi^p E}VV    @VV\varphi_EV \\
   \strc @> id >> \strc \\
   \end{CD} $$ \\
   \end{prop}

This result explains why the big Chern classes commute with
  Adams operations [8] without recourse to hands on computation as
  was done in [8]. \\

\begin{proof}

Let $\hkr^k$ denote the composite of $\hkr$ with the inclusion of
$\texal{X}{k}$ into $\oplus_i \texal{X}{i}$ as a direct summand. Let
$U$ be an affine open subscheme of $X$ and let $v_1,..,v_k$ be
sections of $T_X |_U$. If $(k_1,..,k_p)$ is a $p$-tuple of
nonnegative integers such that $\sum_i k_i = k$, define a
$(k_1,...,k_p)$-multi-shuffle to be a permutation $\sigma$ of
$\{1,..,k\} $ such that \\ $\sigma(1) <...<\sigma(k_1)$,
$\sigma(k_1+1) < ...
<\sigma(k_1+k_2)$,..........,$\sigma(k_1+...+k_{p-1}+1)<...<\sigma(k)$.
Then
  $$ \Delta^p(v_1 \otimes ... \otimes v_k) = \sum_{\{(k_1,....,k_p) | \sum_i k_i =k\} }  \sum_{\sigma \text{ a } (k_1,...,k_p)- \text{shuffle}}
  $$ $$
  \text{sgn}(\sigma)
  v_{\sigma(1)} \otimes .... \otimes v_{\sigma(k_1)} \bigotimes ... \bigotimes v_{\sigma(k_1+...+k_{p-1}+1)} \otimes .. \otimes v_{\sigma(k)} $$
We follow the convention that if $k_i = 0$ then $v_{\sigma(k_1+...+k_{i-1}+1)} \otimes ... \otimes v_{\sigma(k_1+...+k_i)} = 1$. \\

It follows that if $m_p$ denotes the $p$-fold multiplication on $\dd^{\bullet}(X)$ then ,\\

$$m_p \circ \Delta^p (v_1 \otimes ... \otimes v_k) = \sum_{\{(k_1,....,k_p) | \sum_i k_i =k \}} \text{   } \sum_{\sigma \text{ a } (k_1,...,k_p)- \text{shuffle}}
$$ $$  \text{sgn}(\sigma) v_{\sigma(1)} \otimes ... \otimes v_{\sigma(k)} $$

Thus,  $$ m_p \circ \Delta^p (\sum_{\tau \in S_k} \text{sgn}(\tau)v_{\tau(1)} \otimes ... \otimes v_{\tau(k)}) =  $$ $$\sum_{\tau \in S_k}
\text{   } \sum_{\{ (k_1,....,k_p) | \sum_i k_i =k \} } \text{   } \sum_{\sigma \text{ a } (k_1,...,k_p)- \text{shuffle}} $$ $$ \text{sgn}(\tau)
\text{sgn}(\sigma) v_{\sigma(\tau(1))} \otimes ... \otimes v_{\sigma(\tau(k))} $$ Since there are $p^k$  $(k_1,..,k_p)$-shuffles such that
$\sum_i k_i = k$,

$$\sum_{\tau \in S_k} \text{   }
\sum_{\{(k_1,....,k_p) | \sum_i k_i =k \} } \text{     } \sum_{\sigma \text{ a } (k_1,...,k_p)- \text{shuffle}} $$ $$ \text{sgn}(\tau)
\text{sgn}(\sigma) v_{\sigma(\tau(1))} \otimes ... \otimes v_{\sigma(\tau(k))}  = $$ $$ p^k \sum_{\tau \in S_k} \text{sgn}(\tau)v_{\tau(1)}
\otimes ... \otimes v_{\tau(k)} $$ .

From the fact that $\hkr(v_1 \wedge... \wedge v_k) = \frac{1}{k!} \sum_{\tau \in S_k} \text{sgn}(\tau)v_{\tau(1)} \otimes ... \otimes
v_{\tau(k)}$ and the fact that $m_p \circ \Delta^p = \psi^p$, it follows that $\psi^p \circ \hkr^k = p^k \hkr^k$. \\

Now, $\varphi_E \circ \hkr^k= \chn_k(E)$ by Theorem 4,part (iv)
.Thus \\ $\varphi_E \circ \psi^p \circ \hkr^k = p^k \chn_k(E) =
\chn_k(\psi^p E) = \varphi_{\psi^p E} \circ \hkr^k$. This together
with the facts that $\oplus_k \hkr^k = \hkr$ and that $\hkr$ is a
quasi-isomorphism and
therefore, an isomorphism in $\dcat$ prove the desired proposition. \\

\end{proof}

The following corollary is now immediate \\

  \begin{cor}
 The following diagram commutes \\
 $$ \begin{CD}
 \oplus_k T_X^{\otimes k} [-k] @ > \psi^p >> T_X^{\otimes k} [-k]
 \\
 @VV{\oplus_k \ttt_k(\psi^p E)}V     @V{\oplus_k \ttt_k(E)}VV \\
 \strc @> id >> \strc \\
 \end{CD} $$ \\
 \end{cor}

 This is the statement that the big Chern classes commute with
 Adams operations. \\

\section{A formula for the big Chern classes}

This section extends proves a formula for the big Chern classes in terms of the components of the Chern character for vector bundles over an
arbitrary smooth scheme over a field of characteristic $0$. The existence of such a formula was proven in my thesis [8] for smooth
projective varieties using the existence of an ample line bundle together with combinatorial arguments. The method used here is very different from
that of [8]. It is also more general, and works for vector bundles over smooth complex manifolds as well. \\

\subsection{ $\mu \circ \frac{\omega}{1-e^{-\omega}}$ as an element in the PROP $\stdp$}

\subsubsection{A proposition}

We note that there is a PROP $\stdp$ where \\ $\stdp(n,m) :=
\text{Hom}_{\dcat}(T_X^{\otimes n}[-n], T_X^{\otimes m}[-m])$. Let
$\varphi \in \stdp(n,m)$ and \\ $\zeta \in \stdp(n,p)$ . Clearly, we
have a composition $\odot:  \stdp(m,p)  \otimes \stdp(n,m)$ which
takes $\zeta \otimes \varphi$ to $\zeta \circ \varphi$. A
permutation $\sigma$ of $S_n$ gives rise to to elements of
$\stdp(n,n)$ : $l(\sigma)$ is induced at the level of complexes by
the map \\ $v_1 \otimes ... \otimes v_n \leadsto \text{sgn}(\sigma)
v_{\sigma^{-1}(1)} \otimes.... \otimes v_{\sigma^{-1}(n)}$ and
$r(\sigma)$ is induced by \\ $v_1 \otimes ... \otimes v_n \leadsto
\text{sgn}(\sigma) v_{\sigma(1)} \otimes.... \otimes v_{\sigma(n)}$
for sections $v_1,...,v_n$ of $T_X$ over an open affine subscheme
$U$ of $X$ . If $\sigma \in S_n$ and $\tau \in S_m$ then
$\sigma(\varphi) := \varphi \circ l(\sigma)$ and $(\varphi)\tau := r(\tau) \circ \varphi$. \\

We also have a juxtaposition map $$ \times: \stdp(n,m) \otimes
\stdp(n',m') \rar \stdp(n+n',m+m')$$. If $ \eta \in
\stdp(n+n',m+m')$, then $\varphi
* \eta := \varphi \otimes \eta$ \\

Let $\pll$ denote the PROP generated by the Lie operad. \\

We recall from [6](Theorem 3.5.1)  that there is a morphism of PROPs

$$\Upsilon: \pll \rar \stdp $$ so that $$[x1,x2] \in \pll(2,1) \rar \alpha_{T_X} \in \stdp(2,1) := \text{Hom}_{\dcat}(T_X[-1]^{\otimes 2},T_X[-1]) $$

We note that Theorem 2 tells us that the bracket $[,]$ on
$L(\dd^1(X))$ is really \\ $\alpha_{T_X} \in \stdp(2,1)$ after
identifying $L(\dd^1(X))$ with $T_X[-1]$ in $\dcat$ via the map
$\beta$ described in Corollary 3, Section 4 . Let $\omega$ be as in
Theorem 1.We now want to look at $\omega$ as
an element of $\oplus_{1 \leq m \leq n} \stdp(n,m)$. \\

Let $\{z_i\},y,  \{d_i\},d $ be as in the proof of Theorem 1 . Let
\\ $\hat{\omega}: \ufrel{L(\dd^1(X))}  \otimes L(\dd^1(X)) \rar
\ufrel{L(\dd^1(X))} \otimes L(\dd^1(X))$ be the map \\ $z_1 \otimes
.... \otimes z_k \otimes y  \leadsto \sum_i
{(-1)}^{d_i(d_{i+1}+...+d_k)} z_1 \otimes ...
\widehat{z_i}...\otimes z_k \otimes [z_i ,y] $ .  \\ Let
$\hat{\mu}(z_1 \otimes ... \otimes z_k \otimes y)=  \frac{1}{k}
\sum_i {(-1)}^{d(d_{i+1} + .... + d_k)}
z_1 \otimes ... \otimes z_i \otimes y \otimes z_{i+1} \otimes ... \otimes z_k$. \\

Recall the definition of the map \\ $J(L(\dd^1(X))):
\sss^{\bullet}(L(\dd^1(X))) \rar \ufrel{L(\dd^1(X))} \otimes
L(\dd^1(X)) $ from Section 7.1 (Point 4 of pages 26-27). Unlike
Section 7, we will denote $J(L(\dd^1(X)))$ by $B$ to avoid confusing
it with $\bhkr: \oplus_n T_X^{\otimes n}[-n] \rar
\dd^{\bullet}(X)$.    \\

\begin{prop}

The following diagrams commute in $\bbcc$, and hence in $\dcat$ \\

$$ \begin{CD}
\ufrel{L(\dd^1(X))} \otimes L(\dd^1(X)) @> \hat{\omega} >>\ufrel{L(\dd^1(X))} \otimes L(\dd^1(X))\\
 @AA{B \otimes id}A                     @AA{B \otimes \id}A        \\
 \sss^{\bullet}(L(\dd^1(X))) \otimes L(\dd^1(X)) @> \omega >>   \sss^{\bullet}(L(\dd^1(X))) \otimes L(\dd^1(X)) \\
 \end{CD} $$

$$ \begin{CD}
 \ufrel{L(\dd^1(X))} \otimes L(\dd^1(X))  @> \hat{\mu} >> \ufrel{L(\dd^1(X))} \\
                   @AA{B \otimes \id}A        @A{B}AA \\
 \sss^{\bullet}(L(\dd^1(X))) \otimes L(\dd^1(X)) @> \mu >>  \sss^{\bullet}(L(\dd^1(X))) \\
 \end{CD} $$

 \end{prop}

\begin{proof} Since $z_1.....z_k =  z_1... \hat{z_j} ...z_kz_j$ up to a sign, the coefficient of $[z_k,y]$ in $\hat{\omega} \circ (B \otimes id) (z_1... z_k)$ will be
equal to  the coefficient of $[z_j,y]$ in $\hat{\omega} \circ (B \otimes id)$ up to the same sign. The same observation holds with $(B \otimes
id) \circ \omega (z_1...z_k)$ instead of $\hat{\omega} \circ (B \otimes id)$. We thus, need to compare the coefficient of $[z_k,y]$ in
$\hat{\omega} \circ
(B \otimes id) (z_1....z_k)$ and $(B \otimes id) \circ \omega (z_1...z_k)$.\\

In the second map, the coefficient of $[z_k,y]$ is simply
$B(z_1....z_{k-1})$. For a permutation $\sigma \in S_k$ let
$s(\sigma)$ be the sign such that $z_1....z_k = s(\sigma)
z_{\sigma(1)}...z_{\sigma(k)} $. Then \\ $B(z_1....z_k) =
\frac{1}{k!} \sum_{\sigma \in S_k} s(\sigma) z_{\sigma(1)} \otimes
.... \otimes z_{\sigma(k)}$. Let $\tau \in S_{k-1}$ let $\tau_i$ be
the permutation in $S_k$ such that $\tau_i(j) = \tau(j)$ for $j \leq
i-1$, $\tau_i(i) = k$ and $\tau_i(j)= \tau(j-1)$ for $j >i$. Then
$s(\tau_i) = {(-1)}^{d_k(d_{\tau(i+1)+ ... + d_{\tau(k-1)}})}
s_{\tau}$. Thus, the coefficient of $[z_k,y]$ in $\hat{\omega} \circ
(B \otimes id)$ is $$ \frac{1}{k!} \sum_{\tau \in S_{k-1}} \sum_i
{(-1)}^{d_k(d_{\tau(i+1)+ ... + d_{\tau(k-1)}})} s(\tau_i)
z_{\tau(1)} \otimes .... \otimes z_{\tau(k-1)} =  \frac{1}{k!}
\sum_{\tau \in S_{k-1}} \sum_i s(\tau) z_{\tau(1)} \otimes ...
\otimes z_{\tau(k-1)}$$ $$ = \frac{1}{(k-1)!} \sum_{\tau \in
S_{k-1}}  s(\tau) z_{\tau(1)} \otimes .... \otimes z_{\tau(k-1)} $$
This is just $B(z_1...z_{k-1})$.  \\

This shows that the first square commutes. The commuting of the second square is checked by a similar easier calculation \\

\end{proof}

\subsubsection{$\hat{\omega}$ as an element in $\stdp$}

Let $\sigma_{k,n}$ be the permutation of $\{1,...,n\}$ such that $\sigma_{k,n}(n) = k$ and $\sigma_{k,n}(j) = j$ if $j < k$ and $\sigma_{k,n}(j)
= j+1$ otherwise. The map $\sigma \leadsto l(\sigma)$ gives ua a map from $KS_n$ to $\Pi(n,n)$ for any PROP $\Pi$ .\\

The map $\hat{\omega}: \ufrel{L(\dd^1(X))} \otimes L(\dd^1(X)) \rar \ufrel{L(\dd^1(X))} \otimes L(\dd^1(X))$ is clearly the action of the
following element of
$\oplus_{n,m} \stdp(n,m)$ \\

$$ \hat{\omega} : (\sum_{r \geq 1} id \times^{r-1} \times \Upsilon([x_1, x_2])) \circ \sum_n \sum_{k=1}^{k=n} l(\sigma_{k,n}^{-1})  \in \oplus_{1 \leq m < n}
\stdp(n,m) $$ Note that $\hat{\omega}  = \Upsilon(\zeta)$ where $\zeta \in \oplus_{1 \leq m < n} \pll(n,m)$ is given by $$ (\sum_{r \geq 1} id
\times^{r-1} \times [x_1, x_2]) \circ \sum_n \sum_{k=1}^{k=n} l(\sigma_{k,n}^{-1})$$

By convention, if $p \neq n$ then \\ $\circ: \stdp(p,q) \otimes \stdp(m,n) \rar \oplus \stdp(m,n)  = 0$ \\

With this convention, we can think of $\zeta^{\circ k} \in \oplus_{1
\leq m < n} \pll(n,m)$ and \\ $\hat{\omega}^k \in \oplus_{1 \leq m <
n}
\stdp(n,m)$. Denote $\zeta^{\circ k}$ by $\zeta^k$. \\
Clearly, $\hat{\omega}^k = \Upsilon(\zeta^k)$.  This enables us to look at any power series in $\hat{\omega}$ as an element of $\oplus_{1 \leq m
< n} \stdp(n,m)$ . Similarly, any power series in $\zeta$ can be seen as an element in $\oplus_{1 \leq m < n}
\pll(n,m)$\\

We also note that $\hat{\mu} = \oplus_n \frac{1}{n} \sum_k l(\sigma_{k,n}) \in \oplus_m \stdp(m,m) $. \\

Recalling once again that $L(\dd^1(X)) = T_X[-1]$ as Lie algebras in $\dcat$, we see that the map $\hat{\mu} \circ
\frac{\hat{\omega}}{1-e^{-\hat{\omega}}}: \ufrel{L(\dd^1(X))} \otimes L(\dd^1(X)) \rar \ufrel{L(\dd^1(X))}$ is just the action of the element
$$\Theta := \hat{\mu} \circ \frac{\hat{\omega}}{1-e^{-\hat{\omega}}} \in \oplus_{1 \leq m < n} \stdp(n,m)$$ Note that $\Theta = \Upsilon(\hat{\Theta})$ where
$$\hat{\Theta} := \hat{\mu} \circ \frac{\zeta}{1-e^{-\zeta}} \in \oplus_{1 \leq m < n} \pll(n,m) $$.

Let $p:\oplus_i \texal{X}{i} \rar T(T_X[-1])$ be as in Section 8. Since $L(\dd^1(X)) = T_X[-1]$ in $\dcat$, we get the following corollary of proposition 16 \\

\begin{cor}

The following diagram commutes in $\dcat$ \\

$$\begin{CD}
T(T_X[-1]) \otimes T_X[-1] @> \Theta >> T(T_X[-1]) \\
@AA{p \otimes id}A @A{p}AA \\
\oplus_i \texal{X}{i} \otimes T_X[-1] @> \mu \circ \frac{\bar{\omega}}{1-e^{-\bar{\omega}}} >> \oplus_i \texal{X}{i} \\
\end{CD}$$ \\

\end{cor}

\subsection{A return to Theorem 1}

Corollary 1 can be rephrased to say that the following diagram
commutes in $\dcat$.  \\

$$\begin{CD} \dd^{\bullet}(X) \otimes \dd^{\bullet}(X) @> m >> \dd^{\bullet}(X) \\
@AA \hkr \otimes \hkr A            @A \hkr AA \\
\oplus_i \texal{X}{i} \otimes T_X[-1] @> \mu \circ
\frac{\bar{\omega}}{1-e^{-\bar{\omega}}} >> \oplus_i \texal{X}{i} \\
\end{CD}$$ \\

 Let $m_k$ denote the $k$-fold multiplication on $\dd^{\bullet}(X)$. It follows from the above commutative diagram
  that the following diagram commutes $\forall \text{ } k $. \\

$$\begin{CD} \dd^{\bullet}(X) \otimes \dd^{\bullet}(X) @> m >> \dd^{\bullet}(X) \\
@AA{\hkr \otimes (m_k \circ \hkr^{\otimes k})}A            @A{\hkr}AA \\
\oplus_i \texal{X}{i} \otimes T_X^{\otimes k}[-k] @> {(\mu \circ \frac{\bar{\omega}}{1-e^{-\bar{\omega}}} \otimes \id^{k-1}) \circ .... \circ
(\mu \circ
\frac{\bar{\omega}}{1-e^{-\bar{\omega}}}) } >> \oplus_i \texal{X}{i} \\
\end{CD}$$ \\

The following diagram commutes in $\dcat$ since it does so in $\bbcc$ \\

$$\begin{CD}
\dd^{\bullet}(X)  @>{ 1 \otimes (\mu \circ \hkr^{\otimes k})} >> \dd^{\bullet}(X) \otimes \dd^{\bullet}(X) \\
 @AA{m_k \circ \hkr^{\otimes k}}A      @A{\hkr \otimes (m_k \circ \hkr^{\otimes k})}AA \\
T_X^{\otimes k}[-k]     @> {1 \otimes id} >> \oplus_i \texal{X}{i} \otimes T_X^{\otimes k}[-k]\\
\end{CD} $$

Combining this diagram with the previous one we get the following corollary \\

\begin{cor}
The following diagram commutes in $\dcat \text{ } \forall \text{ } k$
\\

$$\begin{CD} \dd^{\bullet}(X) @> id >> \dd^{\bullet}(X) \\
@AA{(m_k \circ \hkr^{\otimes k})}A            @A{\hkr}AA \\
 T_X^{\otimes k}[-k] @> {(\mu \circ \frac{\bar{\omega}}{1-e^{-\bar{\omega}}} \otimes \id^{k-1}) \circ .... \circ (\mu \circ
\frac{\bar{\omega}}{1-e^{-\bar{\omega}}}) \circ (1 \otimes id)}  >> \oplus_i \texal{X}{i} \\
\end{CD}$$
\\

\end{cor}

By Corollary 7 the following diagram commutes in $\dcat$ \\

$$\begin{CD}
T(T_X[-1]) \otimes T_X^{\otimes k}[-k] @> (\Theta \times id^{\times k-1}) \circ .... \circ \Theta >> T(T_X[-1]) \\
@AA{p \otimes id}A @A{p}AA \\
\oplus_i \texal{X}{i} \otimes T_X^{\otimes k}[-k] @> (\mu \circ \frac{\bar{\omega}}{1-e^{-\bar{\omega}}} \otimes \id^{k-1}) \circ .... \circ
(\mu \circ
\frac{\bar{\omega}}{1-e^{-\bar{\omega}}}) \circ (1 \otimes id)  >> \oplus_i \texal{X}{i} \\
\end{CD}$$ \\

Further, the following diagram commutes in $\dcat$ \\

$$ \begin{CD}
T_X^{\otimes k}[-k] @> 1 \otimes id >> T(T_X[-1]) \otimes T_X^{\otimes k}[-k]\\
@AA{id}A  @AA{p \otimes id}A\\
T_X^{\otimes k}[-k]  @> 1 \otimes id >> \oplus_i \texal{X}{i} \otimes T_X^{\otimes k}[-k]\\
\end{CD} $$

Note that the upper morphism in the above diagram just expresses $T_X^{\otimes k}[-k]$ as a summand of $T(T_X[-1]) \otimes T_X^{\otimes k}[-k]$
. We can therefore conclude from the above diagram and the one before that \\

\begin{cor}
The following diagram commutes in $\dcat$ \\

$$\begin{CD}
T_X^{\otimes k}[-k] @> (\Theta \times id^{\times k-1}) \circ .... \circ \Theta >> T(T_X[-1]) \\
@A{id}AA   @A{\gamma}AA \\
T_X^{\otimes k}[-k] @> {(\mu \circ \frac{\bar{\omega}}{1-e^{-\bar{\omega}}} \otimes \id^{k-1}) \circ .... \circ (\mu \circ
\frac{\bar{\omega}}{1-e^{-\bar{\omega}}}) \circ (1 \otimes id)}  >> \oplus_i \texal{X}{i} \\
\end{CD}$$
\\

\end{cor}

Let  $$\Psi_k :=  (\Theta \times id^{\times k-1}) \circ .... \circ \Theta  \in \oplus_{m \leq n} \stdp(n,m) $$  Note that $\Psi_k =
\Upsilon(\hat{\Psi_k})$ where $$\hat{\Psi_k} := (\hat{\Theta} \times id^{\times k-1}) \circ .... \circ \hat{\Theta} \in
\oplus_{m \leq n} \pll(n,m) $$ \\
Let $\Psi_{kl}$ denote the component of $\Psi$ in $\stdp(k,l)$ . Let $\pi: T(T_X[-1]) \rar \oplus_i \texal{X}{i}$ be the standard projection.\\

\begin{thm5}

Let $X$ be a smooth scheme over a field of characteristic $0$. Let
$E$ be a vector bundle on $X$. Then, \\

(i) $\tilde{\ttt_k}(E) = \tilde{\chn_k(E)}\circ \pi + \sum_{l <k}  \tilde{\chn_l(E)} \circ \pi \circ \Psi_{kl} $ \\

(ii) $\ttt_k(E) =  \chn_k(E) \circ \pi + \sum_{l < k} \chn_l(E) \circ \pi \circ \Psi_{kl}$
\\

\end{thm5}

\begin{proof}
Note that  $\pi \circ p = id$. By this observation and by Corollary 9
$$p \circ \Psi = {(\mu \circ \frac{\bar{\omega}}{1-e^{-\bar{\omega}}} \otimes \id^{k-1}) \circ .... \circ (\mu \circ
\frac{\bar{\omega}}{1-e^{-\bar{\omega}}}) \circ (1 \otimes id)} :  T_X^{\otimes k}[-k] \rar  \oplus_i \texal{X}{i}$$ in $\dcat$ \\

Let $\bhkr: \oplus_k T_X^{\otimes k}[-k] \rar \dd^{\bullet}(X)$ and
$\bhkr^k$ be as in Theorem 4. We note that \\ $m_k \circ
\hkr^{\otimes k} = \bhkr^k$. We also recall that for any vector
bundle $E$, we have a morphism $\theta_E: \dd^{\bullet}(X) \rar
\ennnn(E)$ in $\dcat$ so that $\theta_E \circ \bhkr^k =
\tilde{\ttt_k}(E)$ and $\theta_E \circ \hkr = \tilde{\chn}(E)$ by
Theorem 4.
By these observations and Corollary 8,  $$ \tilde{\ttt_k}(E) = \tilde{\chn}(E) \circ \pi \circ \Psi_k $$\\

Note that as $\Psi_k \in \oplus_{m \leq n} \stdp(n,m)$, the contribution of $\chn_l(E) \circ \pi \circ \Psi_{kl}$ to $\ttt_k(E)$ vanishes when
$l>k$ . It only remains to show that $\pi \circ \Psi_{kk} = \pi$. \\

Note that $\Theta = \hat{\mu} \circ \sum_k c_k \hat{\omega}^k$ for
some constants $c_k$. Also note that \\ $\hat{\omega} \in \oplus_{m
< n} \stdp(n,m)$. From the observation that composing an element of
\\ $\oplus_{m \leq n} \stdp(n,m)$ with one of $ \oplus_{m < n}
\stdp(n,m)$ gives an element of \\ $ \oplus_{m < n} \stdp(n,m)$, it
follows that the only component of $\Theta$ in $\oplus_m \stdp(m,m)$
is $\hat{\mu}$. It follows that the component of $\Psi_k$ in
$\oplus_m \stdp(m,m)$ is $(\hat{\mu} \times \id^{\times k-1}) \circ
.... \circ \hat{\mu}$. This map applied to $T_X^{\otimes k}[-k]$ is
just the symmetrization map from $T_X^{\otimes k}[-k]$ to itself. It
follows that the contribution of $\tilde{\chn_k}(E) \circ \pi \circ
\Psi_{kk}$ to
$\tilde{\ttt_k}(E)$ is precisely $\tilde{\chn}_k(E) \circ \pi = \tilde{\chn}_k(E)$. This proves (i) \\

(ii) follows immediately from (i) and from the facts that if $tr: \ennnn(E) \rar \strc$ is the trace map, $\ttt_k(E) = tr \circ
\tilde{\ttt_k}(E)$ and $\chn_k(E) = tr \circ \tilde{\chn}_k(E)$.  \\

\end{proof}

Remark: This theorem gives a formula for the Big Chern classes in terms of the components of the Chern character for arbitrary smooth schemes.
The same proof will go through for complex manifolds as well. This generalizes a similar , more vaguely stated formula $\ttt_k$ in terms of
$\chn_l$ for $l \leq k$ that I obtained in my theses for vector bundles over smooth {\it projective} varieties by some combinatorial methods
[8]. The method there makes it difficult to see the explicit formula for $\Psi_{kl}$. It also requires the existence of an ample line bundle on
the variety for which we are deducing this formula. Even in the smooth projective case, it is difficult to see (i) of Theorem 5 using the methods of [8]. \\

\subsection{Proper subfunctors of the Hodge functors $\caph^q(X,
\Omega^p)$ , $p$,$q \geq 2$}

The formula  for $\ttt_k$ in terms of $\chn_l$ , $ l \leq k$ also easily gives us a method for finding an increasing chain of proper
contravariant subfunctors of the Hodge functors $\caph^q(X, \Omega^p)$ for smooth schemes over a field of characteristic $0$. We note that
 $$\Psi_{kl} \in \text{Hom}_{\dcat}(T_X^{\otimes k}[-k], T_X^{\otimes l}[-l]) = \text{Hom}_{\dcat}(\Omega^{\otimes l}[l], \Omega^{\otimes k}[k])$$
Further, $\pi: T_X^{\otimes l}[-l] \rar \texal{X}{l}$ is identified with $\gamma . k! : \exal{X}{l} \rar \Omega^{\otimes l}[l]$ where $\gamma$
is the symmetrization map . Henceforth , in this subsection, we think of $\chn^l(E)$ and $\ttt_l(E)$ as elements in $\text{Hom}_{\dcat}(\strc,
\exal{X}{l})$ and $\text{Hom}_{\dcat}(\strc,
T_X^{\otimes l}[l])$ respectively. \\

With this convention, the first formula in Theorem 4 can be rewritten to say
$$ \ttt_k(E) = \gamma \circ k! \chn_k(E) + \sum_{l<k} \Psi_{kl} \circ \gamma \circ l! \chn_l(E) $$
The second formula in Theorem 4 may be rewritten in an identical fashion as well, though that does not concern us now. \\

In this picture, $\Psi_{kl} \circ \gamma$ yields a map from
$\caph^l(X ,\Omega^l)$ to $\caph^k(X, \Omega^{\otimes k})$ . Denote
this map by $D_{kl}$. Applying Theorem 4 to $\psi^p E$ , we get \\
$\ttt_k(\psi^p E) = p^k k! \chn_k(E) + \sum_{l<k} p^l
D_{kl}(\chn_l(E)) \text{ } \forall \text{ } p \geq 1$. On the other
hand, \\ $ \ttt_k(\psi^p E) = \sum_{l \leq k} p^l
\ttt_k(\chn^{-1}(\chn_l(E)))$. It follows that $D_{kl}(\chn_l(E)) =
\ttt_k(\chn^{-1}(\chn_l(E)))$ . It was shown in [8] that if $X =
G(r,n)$ a Grassmannian of $r$ dimensional quotient spaces of an $n$
dimensional vector space over a field of characteristic $0$, and if
$n$ is large enough, and if $E = Q$, the canonical quotient bundle
of $X$, then $\ttt_k(\chn^{-1}(\chn_l(E))) \neq 0$ if $l \geq 2$ .
Therefore, the operator $D_{kl}$ does not kill $\caph^{l,l}$ in
general. On the other hand, the Atiyah class $\alpha(T_X) = 0$ is
$X$ is an Abelian
variety (a torus for example). In such a case $D_{kl} = 0$ if $k \neq l$. \\

Therefore, if $X = G(r,n) \times T$ where $T$ is a torus, then $\chn^{-1}(\chn_l(p_1^*Q))$ is not in the kernel of $D_{kl}$. On the other hand,
$D_{kl} (p_2^*Y) =0 \text{ } \forall \text{ } Y \in \caph^{l,l}(T)$. We thus see that $\caph_k^{l,l} := \ker D_{kl}: \caph^{l,l} \rar \caph^k(X
, \Omega^{\otimes k})$ is a proper subfunctor of $\caph^{l,l}$ (as
a theory ) for all $l \geq 2$. \\

Given our current convention,in which we think of $\Psi_{kl}$ as an
element of \\ $\text{Hom}_{\dcat}(\Omega^{\otimes l}[l],
\Omega^{\otimes k}[k])$. Then if $p > q$, $(\Psi_{kq} \otimes
id^{\otimes p-q}) \circ \gamma$  yields a map from $\caph^{p,q}$ to
$\caph^k(X,\Omega^{\otimes k})$. Denote this morphism by $D_{kq}$.
If $p < q$, $D_{kp}$ will denote the map yielded by the element
$\Psi_{kp} \circ \gamma \in \text{Hom}_{\dcat}(\exal{X}{p},
\Omega^{\otimes k}[k])$
from $\caph^{p,q}$ to $\caph^{k-p+q}(X,\Omega^{\otimes k})$.  \\

We see that $\caph_k^{p,q}$ given by $\ker( D_{kq}): \caph^{p,q}
\rar \caph^k(X, \Omega^{\otimes k+p-q})$ if $p > q$ and  \\ $\ker
(D_{kp}): \caph^{p,q} \rar \caph^{k+q-p}(X, \Omega^{\otimes k})$
otherwise is a proper subfunctor (as a theory) of $\caph^{p,q}$. To
see this, again consider the case when $X = G(r,n) \times T$ as
before, $T$ a suitable torus. If $p > q$, and $\alpha_q = \chn^{-1}
\chn_q$, then \\ $D_{kq}
 (\alpha_q(p_1^*Q) \cup p_2^*Y) =( D_{kq}( \chn_q(p_1^*Q)) ) \cup p_2^*Y = \ttt_k(\alpha_qQ) \cup p_2^*Y
\neq 0$ where $Y$ is a nonzero element of $\caph^{p-q,0}(T)$. On the other hand, if $Z \in \caph^{p,q}(T)$ then $D_{kq} ((p_2^*Z)) =0$. This
shows that $\caph_k^{p,q}$ is a proper subfunctor of $\caph^{p,q}$ if $p>q$ . If $p< q$ note that if $Y \in \caph^{0,q-p}(T)$ is nonzero , then
$D_{kp}((\alpha_p(p_1^*Q) \cup p_2^*Y)) = (D_{kp} ( \alpha_p(p_1^*Q))) \cup p_2^*Y = \ttt_k(\alpha_p(p_1^*Q)) \cup p_2^*Y  \neq 0$ and that if
$Z \in \caph^{p,q}(T)$, then $D_{kp}(p_2^*Z) =0$. This proves that $\caph_k^{p,q}$ is a proper subfunctor of
$\caph^{p,q}$ for all $k > q$ where $p,q \geq 2$. \\

\section*{Appendix}

This appendix is meant to collect some facts about graded free Lie algebras used in Lemma 1 and Theorem 1. Proposition 17 is standard. Since I
have not seen Theorem 6 in the literature, I have included it here as a theorem. \\

Let $V$ be a vector space over a field $K$ of characteristic $0$. Let $T(V)$ denote the (graded) tensor algebra generated over $K$ by $V$ in
degree $1$. Let $L(V)$ be the free Lie algebra generated over $K$ by $V$ in degree $1$. \\

Let $I_V: \sss^{\bullet}(L(V)) \rar T(V)$ and $\omega_V: \sss^{\bullet}(L(V)) \otimes L(V) \rar \sss^{\bullet}(L(V)) \otimes L(V)$ be as in
Section 5 (page 16). Then,

\begin{prop} $I_V$ is an isomorphism of graded $K$-vector spaces.\\
\end{prop}

\begin{proof} This is a form of the PBW theorem for $L(V)$ proven in Bahturin [17] . By Theorem 2.10 of Bakhturin [17], $T(V)$ is the universal enveloping algebra of $L(V)$. Let $L(V)_+$ denote the subspace of
$L(V)$ spanned by elements of even degree. Let $L(V)_-$ denote the subspace of $L(V)$ spanned by elements of odd degree. Let $z_1,z_2,.....$ be
a homogenous ordered basis of $L(V)$. We recall from the PBW theorem (Theorem 2.2 of [17]) that the elements $z_{i_i} \otimes ... \otimes
z_{i_n}$ such that $i_j \leq i_{j+1}$ for all $j$ and $z_{i_j} \neq z_{i_{j+1}}$ if $z_{i_j} \in L(V)_-$ form a basis of $T(V)$. Note that the elements
$z_{i_1}....z_{i_n}$ such that $i_j \leq i_{j+1}$ for all $j$ and $z_{i_j} \neq z_{i_{j+1}}$ if $z_{i_j} \in L(V)_-$ form a basis of $\sss^{\bullet}(L(V))$.\\

For a multi-set $S = \{i_1,...,i_n\}$ such that $i_j \leq i_{j+1}$
for all $j$ and $i_j \neq i_{j+1}$ if \\ $z_{i_j} \in L(V)_-$, let
$z_S$ denote the element $z_{i_i} \otimes ... \otimes z_{i_n}$ of
$T(V)$. Let $\pi(z_S)$ denote the element $z_{i_1}...z_{i_n} \in
\sss^{\bullet}(L(V))$. The cardinality $|S|$ of this multi-set is
$n$. Then, let $G_V: T(V) \rar \sss^{\bullet}(L(V))$ be the map such
that $$G_V( \sum_{|S| \leq n, a_S \neq 0 \text{ for some } S \text{
such that } |S|=n } a_S z_S) = \sum_{|S|=n} a_s \pi(z_S)$$  $G_V$ is
a vector space isomorphism by the PBW theorem. Clearly, $G_V \circ
I_V(\pi(z_S)) = \pi(z_S)$ by the proof of Theorem 2.2 of [17]. Since
the $\pi(z_S)$ form a basis of
$\sss^{\bullet}(L(V))$ , $G_V \circ I_V = id$. This proves that $I_V$ is a $K$-vector space isomorphism. \\
\end{proof}

\begin{thm}
The following
diagram commutes in the category of graded $K$ vector spaces. \\

$$\begin{CD}
T(V) \otimes L(V) @> m_V >> T(V) \\
@AA{I_V \otimes id}A @A{I_V}AA \\
\sss^{\bullet}(L(V)) \otimes L(V) @> \mu \circ \frac{\omega_V}{1-e^{-\omega_V}}>> \sss^{\bullet}(L(V)) \\
\end{CD} $$

\end{thm}

\begin{proof}
Let $T(L(V))$ denote the tensor algebra generated over $K$ by
$L(V)$. \\ Let $\hat{\omega}_V:T(L(V))  \bigotimes L(V) \rar T(L(V))
\otimes L(V)$ be the map \\ $z_1 \otimes .. \otimes z_k \bigotimes y
\leadsto \sum_i {(-1)}^{d_i(d_{i+1} +...+ d_k)} z_1 \otimes ..
\widehat{i} .. \otimes  z_k \bigotimes [z_i,y]$ for homogenous
elements $z_1,..,z_k,$ of $L(V)$ of degrees $d_1,...,d_k,d$
respectively. Let $\hat{\mu}: T(L(V)) \bigotimes L(V) \rar T(L(V))$
be the map such that \\ $z_1 \otimes .. \otimes z_{k-1} \bigotimes
z_k \leadsto \frac{1}{k} \sum_{i=1}^{i=k} {(-1)}^{d_k(d_{i} + ..
+d_{k-1})} z_1 \otimes .. \otimes z_{i-1} \otimes z_k \otimes z_i
\otimes... \otimes z_{k-1}$. For a permutation $\sigma$ of $S_k$,
let $s(\sigma)$ be the sign such that $z_1....z_k = s(\sigma)
z_{\sigma(1)}....z_{\sigma(k)}$. Let \\ $J_V: \sss^{\bullet}(L(V))
\rar T(L(V))$ be the
symmetrization map \\ $z_1....z_k \leadsto \sum_{\sigma \in S_k} s(\sigma) z_{\sigma(1)} \otimes .. \otimes z_{\sigma(k)}$. Then the following diagram commutes \\

{\it Step 1 : An analog of Proposition 16} \\

$$\begin{CD}
T(L(V)) \otimes L(V) @>{\hat{\mu} \circ \frac{\hat{\omega}_V}{1-e^{\hat{\omega}_V}}}>> T(L(V)) \\
@A{J_V \otimes id}AA   @AA{J_V}A \\
\sss^{\bullet}(L(V)) @>{\mu \circ \frac{\omega_V}{1-e^{-\omega_V}}}>> \sss^{\bullet}(L(V)) \\
\end{CD} $$

This follows immediately from the fact that the following two diagrams commute . \\

$$ \begin{CD} T(L(V)) \otimes L(V) @>{\hat{\omega}_V}>> T(L(V)) \otimes L(V) \\
@AA{J_V \otimes id}A @A{J_V \otimes id}AA \\
\sss^{\bullet}(L(V)) \otimes L(V) @>\omega_V>> \sss^{\bullet}(L(V)) \otimes L(V)\\
\end{CD} $$

$$\begin{CD}
T(L(V)) \otimes L(V) @>{\hat{\mu}}>> T(L(V)) \\
@AA{J_V \otimes id}A   @A{J_V}AA \\
\sss^{\bullet}(L(V)) \otimes L(V) @>\mu>> \sss^{\bullet}(L(V)) \\
\end{CD} $$

The proof that the above two diagrams commute is word for word identical to that of Proposition 16 (Section 9) with $L(V)$ replacing
$L(\dd^1(X))$ ,$T(L(V))$ replacing \\ $\ufrel{L(\dd^1(X))}$ and $J_V$ replacing $B$. \\

{\it Step 2 : Reduction to a combinatorial question } \\

 The natural inclusion from $L(V)$ to $T(V)$ induces a map of graded algebras \\ $\varphi:
T(L(V)) \rar T(V)$ such that $\varphi
\circ J_V = I_V$. It follows from this that if $m$ denotes the multiplication in $T(L(V))$ we only need to prove the following assertion \\

$$ m \circ (J_V \otimes id) = {\hat{\mu} \circ \frac{\hat{\omega}_V}{1-e^{\hat{\omega}_V}}}$$

Let $z_1,...,z_k , z_{k+1}$ be homogenous elements of $L(V)$ of degrees $d_1,....,d_k$ respectively. Note that $m \circ (J_V \otimes
id)(z_1...z_k \otimes z_{k+1})$ and ${\hat{\mu} \circ \frac{\hat{\omega}_V}{1-e^{\hat{\omega}_V}}}(z_1...z_k \otimes z_{k+1})$ are in the
$K$-span of $\{ z_{\sigma(1)} \otimes ... \otimes z_{\sigma(k+1)} | \sigma \in S_{k+1} \}$. Denote this subspace of $T(L(V))$ by $\calg W$. \\

Note that $S_{k+1}$ has a right action on $\calg W$ such that for a permutation $\tau \in S_{k+1} $, $$\tau(z_{\sigma(1)} \otimes ... \otimes
z_{\sigma(k+1)}) = s(\tau,\sigma) z_{\sigma(\tau(1))} \otimes ..... \otimes  z_{\sigma(\tau(k+1))} $$ where $s(\tau,\sigma)$ is the sign such
that $z_{\sigma(1)}.....z_{\sigma(k+1)} = s(\tau,\sigma) z_{\sigma(\tau(1))} .....  z_{\sigma(\tau(k+1))}$ in \\ $\sss^{k+1}(L(V))$. \\

Let $\sigma(i,l,k+1)$ be the permutation in $S_{k+1}$ such that $\sigma(i,l,k+1)(j)= j$ for $j \leq i-1$. $\sigma(i,l,k+1)(i-1+k) = n-l+k $ for
$1 \leq k \leq l$ and $\sigma(i,l,k+1)(j) = j-l $ for $j \geq i+l$. Let $\nu(i,l,k+1)$ denote the inverse of $\sigma(i,l,k+1)$ in $S_{k+1}$. Let
$\tau_l$ denote the $l$-cycle $(k-l+2 \text{ }  k-l+1 \text{ ... } k+1)$.  \\

{\it Observation 1}

 Identifying $T(L(V)) \otimes L(V)$ as a direct summand of $\oplus_k L(V)^{\otimes k+1}$ of $T(L(V))$, we by a direct
computation that
$$\hat{\mu} \circ \hat{\omega}_V^j = ( \sum_{i=1}^{i=k+1-j} \sigma(i,j+1,k+1)) \circ (id- \tau_{j+1}) \circ ( \sum_{i=1}^{i=k+2-j} \nu(i, j,k+1))  \circ
...$$ $$.. \circ (\sum_{i=1}^{i=k} \sigma(i,2,k+1)) \circ (id- \tau_{2}) \circ  (\sum_{i=1}^{i=k+1} \nu(i, 1,k+1)) $$ on $\calg W$.

Note that $\hat{\mu} \circ  \frac{\hat{\omega}_V}{1-e^{\hat{\omega}_V}} = \sum_j c_j \hat{\mu} \circ \hat{\omega}_V^j$ where $\frac{y}{1-e^{-y}}
= \sum_j c_jy^j$. The above formula thus enables us to express $ \hat{\mu} \circ
\frac{\hat{\omega}_V}{1-e^{\hat{\omega}_V}}$ as the action of an explicit element in the group ring of $S_{k+1}$ on $\calg W$. \\

{\it Observation 2}

On the other hand, $$m \circ (J_V \otimes id)(z_1...z_k \otimes z_{k+1}) = \sum_{\varphi \in S_k} \frac{1}{k!} s(\varphi) z_{\varphi(1)} \otimes
... \otimes z_{\varphi(k)} \otimes z_{k+1}$$

Let $\iota: S_k \rar S_{k+1}$ be the homomorphism fixing $k+1$. It follows from both these observations that we need to prove the following identity in $KS_{k+1}$ for all $k$ \\

$$  \sum_{\varphi \in S_k} \iota(\varphi) = ( \sum_{i=1}^{i=k+1-j} \sigma(i,j+1,k+1)) \circ (id- \tau_{j+1}) \circ ( \sum_{i=1}^{i=k+2-j} \nu(i, j,k+1))  \circ
...$$ $$.. \circ (\sum_{i=1}^{i=k} \sigma(i,2,k+1)) \circ (id- \tau_{2}) \circ  (\sum_{i=1}^{i=k+1} \nu(i, 1,k+1))  \circ \sum_{\varphi \in S_k}
\iota(\varphi) $$

Call this identity (***).

This finishes step 2. \\

{\it Step 3: Proving the combinatorial identity  (***) } \\

Let $W$ (different from $\calg W$) be an infinite dimensional vector
space over $K$ concentrated in degree $0$. Let $L(W)$ and $T(W)$ be
the free Lie algebra generated over $K$ by $W$ and the tensor
algebra generated over $K$ by $W$ respectively. Let
$I_W:\sss^{\bullet}(L(W)) \rar T(W)$ be the symmetrization map such
that $I_W(z_1...z_k)= \sum_{\sigma \in S_k} z_{\sigma(1)} \otimes
... \otimes z_{\sigma(k)}$ for all \\ $z_1,...,z_k \in L(W)$. Let
$\mu$ denote the multiplication on $\sss^{\bullet}(L(W))$ and let
\\ $\omega_W:\sss^{\bullet}(L(W)) \otimes L(W) \rar
\sss^{\bullet}(L(W)) \otimes L(W)  $ be the map such that
$$\omega_W(z_1....z_k \otimes y) = \sum_i z_1....\widehat{i}..z_k
\otimes [z_i,y]$$

It follows from Reutenauer [11] (Chapter 3) that the following diagram commutes in the category of $K$-vector spaces \\

$$\begin{CD}
T(W) \otimes L(W) @>  m >> T(W) \\
@AA{I_W \otimes id}A @A{I_W}AA \\
\sss^{\bullet}(L(W)) \otimes L(W) @> \mu \circ \frac{\omega_W}{1-e^{-\omega_W}}>> \sss^{\bullet}(L(W)) \\
\end{CD} $$

Now, if $z_1,...,z_{k+1}$ are linearly independent elements of $L(W)$, and if $\calg W$ denotes the $K$-span of $\{z_{\sigma(1)} \otimes...
\otimes z_{\sigma(k)} | \sigma \in S_{k+1} \}$ then $S_{k+1}$ has a right action on $\calg W$ such that $\tau(z_{\sigma(1)} \otimes... \otimes
z_{\sigma(k)}) = z_{\sigma(\tau(1))} \otimes... \otimes z_{\sigma(\tau(k))}$. Further, two elements $\alpha$ and $\beta$ in $KS_k$ are equal iff
$\alpha(z_1 \otimes ... \otimes z_n) = \beta(z_1 \otimes .... \otimes z_n)$ in $\calg W$. Let $\iota: S_k \rar S_{k+1}$ be as in the previous
step. The identity (***) of Step 2 follows from this set up once we note that  $$m \circ (I_W \otimes id) (z_1....z_k \otimes z_{k+1}) =
\sum_{\varphi \in S_k} \iota(\varphi)(z_1 \otimes ... \otimes z_{k+1})$$ and $$ I_W \circ (\mu \circ \frac{\omega_W}{1-e^{-\omega_W}})
(z_1....z_k \otimes z_{k+1}) = ( \sum_{i=1}^{i=k+1-j} \sigma(i,j+1,k+1)) \circ (id- \tau_{j+1}) \circ ( \sum_{i=1}^{i=k+2-j} \nu(i, j,k+1))
\circ ...$$ $$.. \circ (\sum_{i=1}^{i=k} \sigma(i,2,k+1)) \circ (id- \tau_{2}) \circ  (\sum_{i=1}^{i=k+1} \nu(i, 1,k+1))  \circ \sum_{\varphi
\in S_k} \iota(\varphi) (z_1 \otimes ... \otimes z_{k+1}) $$

The second of these two identities requires some work. Let $T(L(W))$
be the tensor algebra of $L(W)$. Let \\ $\hat{\omega}_W: T(L(W))
\otimes L(W) \rar T(L(W)) \otimes L(W)$ be the map \\ $z_1 \otimes..
\otimes z_k \bigotimes y \leadsto \sum_i z_1 \otimes ... \widehat{i}
... \otimes z_k \otimes [z_i,y]$. Further, let $\hat{\mu}:T(L(W))
\otimes L(W) \rar T(L(W))$ be the map such that $\hat{\mu}(z_1
\otimes.... \otimes z_{k-1} \bigotimes z_{k}) = \frac{1}{k}
\sum_{i=1}^{i=k} z_1 \otimes .. \otimes z_{i-1} \otimes z_k \otimes
z_i \otimes... \otimes z_{k-1}$. Further, let $J(W):
\sss^{\bullet}(L(W)) \rar T(L(W))$ be the map such that
\\ $z_1...z_k \leadsto \sum_{\sigma \in S_k} z_{\sigma(1)} \otimes.... \otimes z_{\sigma(k)}$.  \\

The following diagrams commute \\
$$ \begin{CD} T(L(W)) \otimes L(W) @>{\hat{\omega}_W}>> T(L(W)) \otimes L(W) \\
@AA{J_W \otimes id}A @A{J_W \otimes id}AA \\
\sss^{\bullet}(L(W)) \otimes L(W) @>\omega_W>> \sss^{\bullet}(L(W)) \otimes L(W)\\
\end{CD} $$

$$\begin{CD}
T(L(W)) \otimes L(W) @>{\hat{\mu}}>> T(L(W)) \\
@AA{J_W \otimes id}A   @A{J_W}AA \\
\sss^{\bullet}(L(W)) \otimes L(W) @>\mu>> \sss^{\bullet}(L(W)) \\
\end{CD} $$

From these, we see that the following diagram commutes \\

$$\begin{CD}
T(L(W)) \otimes L(W) @>{\hat{\mu} \circ \frac{\hat{\omega}_W}{1-e^{\hat{\omega}_W}}}>> T(L(W)) \\
@A{J_W \otimes id}AA   @AA{J_W}A \\
\sss^{\bullet}(L(W)) @>{\mu \circ \frac{\omega_W}{1-e^{-\omega_W}}}>> \sss^{\bullet}(L(W)) \\
\end{CD} $$

Now , a direct computation shows us that \\

$$ ( \hat{\mu} \circ \frac{\hat{\omega}_W}{1-e^{\hat{\omega}_W}}) \circ (J_W \otimes id) (z_1 \otimes .... \otimes z_k \bigotimes z_{k+1})$$ $$ =
( \sum_{i=1}^{i=k+1-j} \sigma(i,j+1,k+1)) \circ (id- \tau_{j+1}) \circ ( \sum_{i=1}^{i=k+2-j} \nu(i, j,k+1)) \circ ...$$ $$.. \circ
(\sum_{i=1}^{i=k} \sigma(i,2,k+1)) \circ (id- \tau_{2}) \circ  (\sum_{i=1}^{i=k+1} \nu(i, 1,k+1))  \circ \sum_{\varphi \in S_k} \iota(\varphi)
(z_1 \otimes ... \otimes z_{k+1}) $$

Thus,
$$ I_W \circ (\mu \circ \frac{\omega_W}{1-e^{-\omega_W}})
(z_1....z_k \otimes z_{k+1}) = ( \sum_{i=1}^{i=k+1-j} \sigma(i,j+1,k+1)) \circ (id- \tau_{j+1}) \circ ( \sum_{i=1}^{i=k+2-j} \nu(i, j,k+1))
\circ ...$$ $$.. \circ (\sum_{i=1}^{i=k} \sigma(i,2,k+1)) \circ (id- \tau_{2}) \circ  (\sum_{i=1}^{i=k+1} \nu(i, 1,k+1))  \circ \sum_{\varphi
\in S_k} \iota(\varphi) (z_1 \otimes ... \otimes z_{k+1}) $$

Thereby proving Theorem 5. \\

\end{proof}

\section*{Acknowledgements}

I thank Prof. Victor Ginzburg for introducing me to the works of
Caldararu [10] and Markarian [3] and asking whether the material in
these works had anything to do with the relation between the big
Chern classes and the Chern character. In addition, Prof. Ginzburg
went through the first draft of this work and gave me many helpful
suggestions. I am grateful Prof. Madhav Nori for many useful
discussions. I also thank Prof. Kaan Akin , Dr. Victor Protsak, and
the referees for the first and current versions of this paper for
their useful suggestions.

\end{document}